\documentclass[a4paper,reqno,12pt]{amsart}
\usepackage{enumitem}
\usepackage{amsmath,amsfonts,amssymb}
\usepackage{graphics,color}
\usepackage[latin1]{inputenc}
\usepackage[colorlinks = true,
            linkcolor = blue,
            urlcolor  = blue,
            citecolor = blue,
            anchorcolor = blue]{hyperref}

\newtheorem{theorem}{Theorem}[section]
\newtheorem{lemma}[theorem]{Lemma}

\makeatletter
\@namedef{subjclassname@2020}{\textup{2020} Mathematics Subject Classification}
\makeatother

\allowdisplaybreaks

\numberwithin{equation}{section}
\newtheorem{proposition}[theorem]{Proposition}
\newtheorem{corollary}[theorem]{Corollary}

\newtheorem{remark}[theorem]{Remark}

\newcommand{\R}{\mathbb{R}}
\newcommand{\dd}{\,dx}
\newcommand{\dist}{\operatorname{dist}}
\newcommand{\spanop}{\operatorname{span}}
\newcommand{\rad}{\operatorname{rad}}

\usepackage{cite}
\numberwithin{equation}{section}

\allowdisplaybreaks[4]
\begin{document}
\title{Optimal Stability for Bivariate Weighted Fourth-Order Inequalities}

\author{Shengbing Deng}
\address{School of Mathematics and Statistics, Southwest University, Chongqing 400715, P. R. China.}
\email{shbdeng@swu.edu.cn}

\author{Bo Zhang$^*$}
\address{School of Mathematics and Statistics, Southwest University, Chongqing 400715, P. R. China.}
\email{zhangbo371013@163.com}

\footnotetext[1]{Corresponding author.}


\begin{abstract}
In this paper, we prove sharp stability estimates for a two-component radial weighted
fourth-order Caffarelli--Kohn--Nirenberg inequality in the two radial regimes for which
the scalar optimizers are explicitly known and the corresponding spectral
non-degeneracy holds. These regimes are
\(4-N<\alpha<2,\ \ell=-\alpha\), and
\(2<\alpha<N,\ \ell=4(\alpha-2)(N-2)/(N-\alpha)-\alpha\). Let $N\ge 3$, $p=2(N+\ell)/(N-4+\alpha)>2$, \(\mu_1,\mu_2\geq0\), \(\kappa>0\), and $r,s>1$ satisfy $r+s=p$. For radial pairs $(u,v)\in X_{\mathrm{rad}}\times X_{\mathrm{rad}}$, we determine the sharp constant $S_{r,s,\mu_1,\mu_2,\kappa}$ in the coupled inequality
\[
S_{r,s,\mu_1,\mu_2,\kappa}
\left(\int_{\mathbb R^N}|x|^\ell\big(\mu_1|u|^p+\mu_2|v|^p+p\kappa |u|^r|v|^s\big)dx\right)^{2/p}
\le
\int_{\mathbb R^N}|x|^\alpha\big(|\Delta u|^2+|\Delta v|^2\big)dx .
\]
The sharp constant is expressed through the scalar sharp constant and the minimum of an explicit one-dimensional ratio function. We also characterize all extremals and prove sharp quantitative stability estimates. The stability exponent is determined by the degeneracy of the ratio minimizers: it is quadratic in the non-degenerate case and becomes square-root type in the degenerate case.
\end{abstract}

\maketitle
\vspace {-1cm}
\section{Introduction}
Given a sharp functional inequality, one may ask whether its deficit controls,
in a quantitative way, the distance of a function to the manifold of optimizers.
This question goes back to Brezis and Lieb~\cite{BL1985} and received a
seminal answer in the work of Bianchi and Egnell~\cite{BE1991} on the
Sobolev inequality. Since then, stability estimates for Sobolev-type inequalities
have been developed in several directions, including remainder terms, gradient
stability, sharp constants, and stability near critical points
\cite{LW2000,FN2019,FZ2022,K2023,K2025,LZ2025,CKW2024,ZW2024,CS2025}.
These works show that the stability theory of sharp Sobolev inequalities is not
merely a refinement of the inequality itself, but also reflects subtle information
on the optimizer manifold, spectral non-degeneracy, and possible loss of
compactness.

Weighted problems originate in the
Caffarelli--Kohn--Nirenberg inequality \cite{CKN1984} and may exhibit symmetry
breaking \cite{FS2003}. Their stability theory now includes
\cite{F2022,CFL2024,WW2022,WW2024,FP2024,DGK2025,ZW2024,WW2026,ZZ2026}. The
fourth-order weighted equations relevant here were analyzed in
\cite{DGT2023,DT2023,DT2024,DT2025}; their radial bubbles and spectral
information provide the scalar input below.

Let \(N\geq 3\). We consider the scalar weighted fourth-order  Caffarelli--Kohn--Nirenberg inequality
\begin{equation}\label{eq001}
 \int_{\R^N}|x|^\alpha |\Delta u|^2 \dd
 \geq
 S_{0}
 \left(
 \int_{\R^N}|x|^\ell |u|^p \dd
 \right)^{2/p},
 \quad u\in X_{\mathrm{rad}},
\end{equation}
where
\[
 p:=\frac{2(N+\ell)}{N-4+\alpha},\quad X_{\mathrm{rad}}:=\overline{C^{\infty}_{c,\mathrm{rad}}\bigl(\mathbb{R}^N \setminus \{0\}\bigr)}^{\| \cdot \|_{\alpha}}
\]
with the corresponding norm
\[
\|u\|^2_{\alpha}:=\int_{\R^N}|x|^\alpha |\Delta u|^2 \dd,
\]
and
\[
S_0:=\inf_{u\in X_{\mathrm{rad}} \backslash \{ 0\} }\frac{\int_{\mathbb{R}^{N}}|x|^{\alpha }|\Delta u|^{2}dx}{\left( \int_{\mathbb{R}^{N}}|x|^{\ell}|u|^{p}dx\right) ^{\frac{2}{p}}}.
\]
We focus on the following two cases, in which the scalar radial optimizers are known explicitly.

\noindent\textbf{Case A.}
\begin{equation}\label{case-A}
 4-N<\alpha<2,\quad \ell=-\alpha,\quad
 p=\frac{2(N-\alpha)}{N-4+\alpha}.
\end{equation}
The positive radial optimizer is
\begin{equation}\label{bubble-A}
 W_\lambda(x)
 =
 C_{N,\alpha}\lambda^{(N-4+\alpha)/2}
 \left(1+\lambda^{2-\alpha}|x|^{2-\alpha}\right)^{-(N-4+\alpha)/(2-\alpha)},
 \quad \lambda>0.
\end{equation}
\noindent\textbf{Case B.}
\begin{equation}\label{case-B}
 2<\alpha<N,\quad
 \ell=\frac{4(\alpha-2)(N-2)}{N-\alpha}-\alpha,\quad
 p=\frac{2(N-4+\alpha)}{N-\alpha}.
\end{equation}
The positive radial optimizer is
\begin{equation}\label{bubble-B}
 W_\lambda(x)
 =
 C_{N,\alpha}\lambda^{(N-\alpha)/2}|x|^{2-\alpha}
 \left(1+\lambda^{\alpha-2}|x|^{\alpha-2}\right)^{-(N-\alpha)/(\alpha-2)},
 \quad \lambda>0.
\end{equation}
In both cases, the family \(W_\lambda\) satisfies the natural scaling relation
\begin{equation}
 W_\lambda(x)=\lambda^{(N-4+\alpha)/2}W_1(\lambda x),
 \label{scaling}
\end{equation}
with the understanding that the explicit expression of \(W_1\) is given by either \eqref{bubble-A} or \eqref{bubble-B}.
Define
\[\mathcal{M}_0:= \{\, c W_{\lambda} : c \in \mathbb{R},\ \lambda > 0 \}.\]
The nonzero minimizers of \eqref{eq001} are precisely $\mathcal{M}_0\backslash\{0\}$.
The scalar input is recorded in a form that includes the compactness statement
used later.

\begin{theorem}\cite{DGT2023,DT2023}\label{thm:scalar-stability}
Let \(N\geq 3\). Assume that one of \eqref{case-A} and \eqref{case-B} holds. Then there exists a constant $B=B(N,\alpha)>0$
such that, for every $u\in X_{\mathrm{rad}}$,
\begin{equation}\label{scalar-stability}
\int_{\mathbb{R}^N}
|x|^\alpha |\Delta u|^2\,dx
-
S_{0}
\left(
\int_{\mathbb{R}^N}
|x|^\ell |u|^p\,dx
\right)^{\frac{2}{p}}
\geq
B\,\operatorname{dist}
\bigl(u,\mathcal{M}_0\bigr)^2.
\end{equation}
where
\[
\operatorname{dist}
\bigl(u,\mathcal{M}_0\bigr)
:=\inf_{c\in\mathbb{R},\,\lambda>0}
\|u-cW_{\lambda}\|_{\alpha}.
\]
\end{theorem}

Our motivation comes from bivariate sharp inequalities. Zhang et al. \cite{ZZZ2026} studied a two-component system Hardy--Sobolev inequality and revealed that the stability exponent may become weaker when the finite-dimensional parameter determining the ratio of the two components is degenerate.

Our first contribution is to determine the sharp constant of the system inequality and characterize its minimizers.

For \(\mu_1,\mu_2\geq0\), \(\kappa>0\), \(r,s>1\), and $ r+s=p$,
we define the two-component nonlinear functional
\begin{equation}
 P(u,v)
 :=
 \int_{\R^N}|x|^\ell
 \left(
\mu_1 |u|^p+\mu_2 |v|^p+p\kappa |u|^r|v|^s
 \right)\dd .
 \label{P-def}
\end{equation}
We use the product
space
\[
\mathcal X:=X_{\mathrm{rad}}\times X_{\mathrm{rad}}
\]
equipped with the norm
\[
\|(u,v)\|_{\mathcal X}^{2}
:=
\|u\|_{\alpha}^{2}+\|v\|_{\alpha}^{2}.
\]
The associated system inequality is written as
\begin{equation} \label{system-ineq}
 \|(u,v)\|_{\mathcal X}^{2}
 \geq
 S_{r,s,\mu_1,\mu_2,\kappa}
 \left(P(u,v)\right)^{2/p},
 \quad (u,v)\in \mathcal X,
\end{equation}
where
\[
S_{r,s,\mu_1,\mu_2,\kappa}:=\inf_{\substack{(u,v)\in\mathcal X \backslash \{ (0,0)\}\\P(u,v)>0 }}\frac{\|(u,v)\|_{\mathcal X}^{2}}{\left(P(u,v)\right)^{2/p}}.
\]

\begin{theorem}\label{th1.2}
Let \(N\geq 3\). Assume that one of \eqref{case-A} and \eqref{case-B} holds. Then, for any
$(u,v)\in \mathcal X$, we have \eqref{system-ineq}. The sharp constant
$S_{r,s,\mu_1,\mu_2,\kappa}$ is given by
\[
 S_{r,s,\mu_1,\mu_2,\kappa}=m_*S_{0},
\]
where
\[
 m_*:=\min_{t\in[0,\infty]}g(t),\quad
g(t):=
 \frac{1+t^2}
 {\left(\mu_1+\mu_2t^p+p\kappa t^s\right)^{2/p}},
\quad t\in[0,\infty),
\]
and the endpoint values are understood in the extended sense:
\[
g(0):=
\begin{cases}
\mu_1^{-2/p}, & \mu_1>0,\\
+\infty, & \mu_1=0,
\end{cases}
\]
\[
g(\infty):=
\begin{cases}
\mu_2^{-2/p}, & \mu_2>0,\\
+\infty, & \mu_2=0.
\end{cases}
\]

Define
\[
 \mathcal S:=\{t\in[0,\infty]:g(t)=m_*\}.
\]
Then the set of minimizers of \eqref{system-ineq} is
\[
\mathcal M
:=
\Big\{
A\big(W_\lambda,\sigma t W_\lambda\big):
A\in\mathbb R\setminus\{0\},\ \lambda>0,\ t\in \mathcal S\cap [0,\infty),\
\sigma\in\{-1,1\}
\Big\}.
\]
If \(\infty\in\mathcal S\), one has in addition the semi-trivial branch
\[
\{(0,AW_\lambda):A\in\mathbb R\setminus\{0\},\lambda>0\}.
\]
In particular, if
\((r,s,\mu_1,\mu_2)=(2,2,2\kappa,2\kappa)\), then
\(\mathcal S=[0,\infty]\), and the above description is equivalently
\begin{equation}\label{eq007}
\mathcal M_{\mathrm{flat}}
:=
\{(aW_\lambda,bW_\lambda):(a,b)\in\mathbb R^2\setminus\{(0,0)\},\lambda>0\}.
\end{equation}
\end{theorem}

Define the deficit functions \(\delta(u)\) and \(\delta(u,v)\) by
\begin{equation}\label{eq004}
\delta(u):=\int_{\mathbb{R}^N}
|x|^\alpha |\Delta u|^2\,dx
-
S_{0}
\left(
\int_{\mathbb{R}^N}
|x|^\ell |u|^p\,dx
\right)^{\frac{2}{p}},
\quad u\in X_{\rad},
\end{equation}
and
\begin{equation}\label{eq005}
\delta(u,v):=\|(u,v)\|_{\mathcal X}^{2}
-S_{r,s,\mu_1,\mu_2,\kappa}
 \left(P(u,v)\right)^{2/p},
 \quad (u,v)\in \mathcal X.
\end{equation}
The degenerate ratio set is
\[
\mathcal S_{\rm deg}
:=
\mathcal S\cap
\Big(
\{t\in(0,\infty):g''(t)=0\}
\cup
\{0:g''_+(0)=0\}
\cup
\{\infty:\widetilde g''_+(0)=0\}
\Big).
\]
The behavior of \(g\) at infinity is described by
\begin{equation}\label{eq006}
\widetilde g(\tau):=g(1/\tau)
=
\frac{1+\tau^2}
{\left(\mu_2+\mu_1\tau^p+p\kappa\tau^r\right)^{2/p}},
\quad \tau\ge0.
\end{equation}
The distance to \(\mathcal M\) is defined by
\[
\operatorname{dist}\big((u,v),\mathcal M\big)
:=
\inf_{(U,V)\in\mathcal M}
\left(
\|u-U\|_\alpha^2+\|v-V\|_\alpha^2
\right)^{1/2}.
\]

Our main result is as follows.

\begin{theorem}\label{main-thm}
Assume that one of \eqref{case-A} and \eqref{case-B} holds.
\begin{itemize}
  \item If $(r,s,\mu_1,\mu_2)\neq(2,2,2\kappa,2\kappa)$,  then there exists a constant
$
 C=C(N,\alpha,\mu_1,\mu_2,\kappa,r,s)>0
$
such that for all nonzero \((u,v)\in \mathcal X\),
\begin{equation}\label{main-stability}
 \left(\frac{\delta(u,v)}{\|(u,v)\|_{\mathcal X}^{2}}\right)^{\iota}
 \geq
 C\,\frac{\dist((u,v),{\mathcal M})^2}{\|(u,v)\|_{\mathcal X}^{2}},
\end{equation}
where
\[
\iota=
\begin{cases}
1, & \mathcal S_{\rm deg}=\emptyset,\\[1mm]
\dfrac12, & \mathcal S_{\rm deg}\neq\emptyset.
\end{cases}
\]
  \item If \(
(r,s,\mu_1,\mu_2)=(2,2,2\kappa,2\kappa),
\)
then there exists a constant
\(
C=C(N,\alpha,\kappa)>0
\)
such that, for every nonzero \((u,v)\in\mathcal X\),
\begin{equation}\label{flat-stability}
\delta(u,v)
\geq
C\,\operatorname{dist}
\bigl((u,v),\mathcal M_{\mathrm{flat}}\bigr)^2.
\end{equation}
\end{itemize}
Moreover, these exponents are optimal in the corresponding cases.
\end{theorem}

\begin{remark}
\label{rem:flat-case}
The exceptional flat case
\(
(r,s,\mu_1,\mu_2)=(2,2,2\kappa,2\kappa)
\)
is qualitatively different from the isolated degenerate cases covered
by Theorem~\ref{main-thm}. Indeed, since the ratio function satisfies
\(
g(t)\equiv (2\kappa)^{-1/2},
\, t\in[0,\infty].
\)
Thus every component ratio is minimizing, and the minimizing manifold is
\(
\mathcal M_{\mathrm{flat}}.
\)
The ratio-changing direction is therefore tangent to the minimizing
manifold rather than a degenerate normal direction. Once this additional
tangent direction is removed, the second variation is coercive on the
normal space, which yields the quadratic estimate
\eqref{flat-stability}.
\end{remark}

The paper is organized as follows. In Section 2, we collect the scalar weighted fourth-order facts and analyze the finite-dimensional function $g$ which determines the ratio of the two components. In Section 3, we prove the sharp system inequality and characterize all extremals. Section 4 is devoted to compactness of minimizing sequences and to a qualitative stability result near the minimizing manifold. In Section 5, we carry out the local stability analysis and prove the sharp estimates in both the non-degenerate and degenerate cases. We also treat the exceptional flat case, identifies the additional
\(O(2)\)-symmetry of the optimizer manifold, and proves the corresponding quadratic stability estimate after projecting out the extra ratio direction. Throughout the paper, \(A\asymp B\) means that there exist constants \(c,C>0\), independent of the local parameters, such that
\(cB\le A\le CB\).

\section{Preliminaries}\label{sec:preliminaries}

We collect the scalar facts used in the proof. Throughout the paper, \(W_\lambda\) denotes the scalar optimizer in either Case A or Case B.
For the scalar radial inequality \eqref{eq001}, the best constant \(S_{0}\) is achieved precisely by
\[
 cW_\lambda,\quad c\in\R\setminus\{0\},\quad \lambda>0.
\]
Moreover,
\begin{equation}
 \|W_\lambda\|_\alpha^2
 =
 S_{0}
 \left(
 \int_{\R^N}|x|^\ell W_\lambda^p\dd
 \right)^{2/p}.
\end{equation}
After a harmless normalization, we assume
\begin{equation} \label{scalar-normalization}
 \int_{\R^N}|x|^\ell W_\lambda^p\dd=1,
 \quad
 \|W_\lambda\|_\alpha^2=S_{0}.
\end{equation}

\begin{proposition}\cite{DT2023,DGT2023}
\label{prop:scalar-gap-A}
Assume that one of \eqref{case-A} and \eqref{case-B} holds. Consider the radial eigenvalue problem
\begin{equation}\label{scalar-eigen-A}
 \Delta(|x|^\alpha\Delta\phi)
 =
\Lambda |x|^\ell W_\lambda^{p-2}\phi,
 \qquad \phi\in X_{\rad}.
\end{equation}
Then the first two radial eigenspaces are generated by \(W_\lambda\) and
\(\partial_\lambda W_\lambda\), respectively, and
\[
 \Lambda_1=S_{0},\qquad
 \Lambda_2=(p-1)S_{0}.
\]
Moreover, there is a radial spectral gap above the scaling direction: if
\[
 \langle \phi,W_\lambda\rangle_\alpha=0,
 \quad
 \langle \phi,\partial_\lambda W_\lambda\rangle_\alpha=0,
\]
then
\begin{equation}\label{scalar-gap-ineq-A}
 \Lambda_3
 \int_{\R^N}|x|^\ell W_\lambda^{p-2}\phi^2\dd
 \leq
 \|\phi\|_\alpha^2
\end{equation}
for some \(\Lambda_3>(p-1)S_{0}\). Here \(\Lambda_3\) denotes the first radial eigenvalue above \((p-1)S_{0}\).
\end{proposition}

We now prove the scalar compactness of minimizing sequences.

\begin{lemma}\label{lem:scalar-cylinder-compactness}
Assume that one of \eqref{case-A} and \eqref{case-B} holds, then every normalized minimizing sequence for
\eqref{eq001} is strongly compact modulo dilation and multiplication by a
nonzero constant.
\end{lemma}

\begin{proof}
Put
\[
 \beta=\frac{N-4+\alpha}{2},\qquad
 c_{\alpha,N}=\frac{(N-4+\alpha)(N-\alpha)}4,
\qquad u(r)=r^{-\beta}w(\log r).
\]
A direct radial computation, using \(\beta p=N+\ell\), gives
\[
\begin{split}
 \|u\|_\alpha^2
 &=|\mathbb S^{N-1}|\int_{\mathbb R}
 \{(w'')^2+[(2-\alpha)^2+2c_{\alpha,N}](w')^2
             +c_{\alpha,N}^2w^2\}\,dt,\\
 \int_{\mathbb R^N}|x|^\ell|u|^p\,dx
 &=|\mathbb S^{N-1}|\int_{\mathbb R}|w|^p\,dt .
\end{split}
\]
The first expression is an equivalent \(H^2(\mathbb R)\)-norm, denoted by
\(\|\cdot\|_{\mathcal H}\), and \(W_\lambda\) becomes a translate of
\(W_1\).

Let \(\|w_n\|_{\mathcal H}=1\) be minimizing. Its \(L^p\)-norm stays away
from zero. For \(I_k=(k,k+1)\) and \(I_k^*=(k-1,k+2)\), local Sobolev
embedding and finite overlap give
\[
 \int_{\mathbb R}|w|^p\,dt
 \le C\left(\sup_k\int_{I_k}|w|^p\,dt\right)^{1-2/p}
       \|w\|_{\mathcal H}^2 .
\]
Indeed, factor out the local mass to the power \(1-2/p\), use the embedding
\(H^2(I_k^*)\hookrightarrow L^p(I_k)\), and sum. Hence, after integer
translations, a subsequence satisfies \(w_n\rightharpoonup w\) in
\(\mathcal H\) and \(w_n\to w\) in \(L^p_{\rm loc}\), with \(w\ne0\).

Let \(\widehat S\) be the transformed sharp constant and put
\(z_n=w_n-w\). Hilbert orthogonality and the Brezis--Lieb lemma yield
\[
\begin{split}
 1&=\|w\|_{\mathcal H}^2+\lim_n\|z_n\|_{\mathcal H}^2,\\
 \widehat S^{-p/2}
 &=\|w\|_{L^p}^p+\lim_n\|z_n\|_{L^p}^p .
\end{split}
\]
The sharp inequality applied to both terms, followed by the strict
concavity of \(q\mapsto q^{2/p}\), forces one of the two \(L^p\)-masses to
vanish. Since \(w\ne0\), the remainder vanishes and \(z_n\to0\) in
\(\mathcal H\). Thus \(w\) is an optimizer. Transforming back and using the
optimizer classification in \cite{DGT2023,DT2023} proves the lemma.
\end{proof}

\begin{lemma}\label{lem:weak-bubble}
For every \(f\in X_{\rm rad}\), one has
\[
\langle f,W_\lambda\rangle_\alpha\longrightarrow 0
\quad\text{as }\lambda\to0\text{ or }\lambda\to\infty.
\]
Equivalently, \(W_\lambda\rightharpoonup0\) weakly in \(X_{\rm rad}\) at both ends of the scale parameter.
\end{lemma}

\begin{proof}
It is enough to prove the assertion first for
\(f\in C^\infty_{c,{\rm rad}}(\mathbb R^N\setminus\{0\})\). Then
\(\Delta f\) is supported in an annulus \(R_1<|x|<R_2\). By the explicit formulae
\eqref{bubble-A} and \eqref{bubble-B}, together with the scaling relation
\eqref{scaling}, the functions \(\Delta W_\lambda\) and their weighted energy
mass either concentrate at a scale disjoint from this fixed annulus or are seen
there only through the decaying tail of the bubble. Consequently,
\[
\int_{R_1<|x|<R_2}|x|^\alpha \Delta f\,\Delta W_\lambda\,dx\to0
\quad\text{as }\lambda\to0\text{ or }\lambda\to\infty.
\]
This proves \(\langle f,W_\lambda\rangle_\alpha\to0\) for such test functions.
For a general \(f\in X_{\rm rad}\), choose
\(f_j\in C^\infty_{c,{\rm rad}}(\mathbb R^N\setminus\{0\})\) such that
\(f_j\to f\) in \(X_{\rm rad}\). This density property is the usual definition
of the homogeneous radial space in the present weighted setting. Since
\(\|W_\lambda\|_\alpha=\|W_1\|_\alpha\),
\[
|\langle f-f_j,W_\lambda\rangle_\alpha|
\le \|f-f_j\|_\alpha\|W_1\|_\alpha .
\]
Letting first \(\lambda\to0\) or \(\lambda\to\infty\), and then \(j\to\infty\),
gives the claim.
\end{proof}
\begin{proof}[Proof of Theorem~\ref{thm:scalar-stability}]
Lemma~\ref{lem:scalar-cylinder-compactness} shows that if
\(\|u_n\|_\alpha=1\) and the scalar deficit tends to zero, there are
\(c_n\ne0\) and \(\lambda_n>0\) such that
\begin{equation}
 u_n-c_nW_{\lambda_n}\longrightarrow0
 \quad\hbox{in }X_{\mathrm{rad}}.
 \label{eq:scalar-compactness-precise}
\end{equation}
This proves the compactness assertion.

We prove the quantitative assertion. Assume first that \(u\) is sufficiently
close to \(\mathcal M_0\) and choose a nearest point
\[
 \dist(u,\mathcal M_0)=\|u-cW_\lambda\|_\alpha,
 \qquad u=cW_\lambda+\varphi .
\]
The pair \((c,\lambda)\) exists by the weak vanishing of \(W_\lambda\) at
both ends of the scale parameter, proved in Lemma~\ref{lem:weak-bubble}; in a
fixed sufficiently small neighborhood one has \(c\ne0\). Differentiating the
squared distance with respect to \(c\) and \(\tau=\log\lambda\) gives
\[
 \langle\varphi,W_\lambda\rangle_\alpha=0,
 \qquad
 \langle\varphi,\lambda\partial_\lambda W_\lambda\rangle_\alpha=0.
\]
By two-homogeneity we divide by \(|c|\), and by dilation invariance we reduce
to \(|c|=1\), \(\lambda=1\). The map
\[
 z\longmapsto
 \left(\int_{\mathbb R^N}|x|^\ell|z|^p\,dx\right)^{2/p}
\]
is twice differentiable at the positive bubble on the weighted \(L^p\)
space. Hence, uniformly in \(c\ne0\) and \(\lambda>0\),
\[
 \delta(cW_\lambda+\varphi)
 =Q_\lambda(\varphi)+o(\|\varphi\|_\alpha^2),
\]
where
\[
 Q_\lambda(\varphi)
 =\|\varphi\|_\alpha^2-(p-1)S_0
 \int_{\mathbb R^N}|x|^\ell W_\lambda^{p-2}\varphi^2\,dx .
\]
By Propositions~\ref{prop:scalar-gap-A},
\[
 Q_\lambda(\varphi)
 \ge
 \left(1-\frac{(p-1)S_0}{\Lambda_3}\right)
 \|\varphi\|_\alpha^2.
\]
The remainder divided by \(\|\varphi\|_\alpha^2\) tends to zero when
\(\|\varphi\|_\alpha/|c|\to0\), independently of \(\lambda\). Thus, in a
fixed conical neighborhood of \(\mathcal M_0\),
\[
 \delta(u)\ge c\dist(u,\mathcal M_0)^2 .
\]
It remains to pass from the local estimate to the global one. If the stated
estimate failed, two-homogeneity would give a sequence with
\(\|u_n\|_\alpha=1\) such that
\[
 \delta(u_n)\to0,
 \qquad
 \frac{\delta(u_n)}{\dist(u_n,\mathcal M_0)^2}\to0.
\]
By \eqref{eq:scalar-compactness-precise}, the distance to \(\mathcal M_0\) tends
to zero. Hence \(u_n\) lies in the uniform local regime for large \(n\), which
contradicts the second limit. This proves \eqref{scalar-stability}.
\end{proof}
We recall the following finite-dimensional classification from \cite{ZZZ2026}.
\begin{lemma}\label{3.0}\cite{ZZZ2026}
Assume
$
(r,s,\mu_1,\mu_2)\neq(2,2,2\kappa,2\kappa).
$
Then \(\mathcal S\) is finite. Moreover, every minimizer \(t_0\in \mathcal S\) is either nondegenerate
or degenerate in the following sense.

\begin{enumerate}
    \item If \(t_0\in[0,\infty)\) is nondegenerate, then
    \[
    g'(t_0)=0,
    \quad
    g''(t_0)>0,
    \]
    where derivatives at \(t_0=0\) are understood as right derivatives.

    If \(t_0=\infty\) is nondegenerate, then
    \[
    \widetilde g'(0)=0,
    \quad
    \widetilde g''(0)>0.
    \]

    \item If \(t_0\in[0,\infty)\) is degenerate, then
    \[
    g'(t_0)=g''(t_0)=g'''(t_0)=0,
    \quad
    g^{(4)}(t_0)>0.
    \]

    If \(t_0=\infty\) is degenerate, then
    \[
    \widetilde g'(0)=\widetilde g''(0)=\widetilde g'''(0)=0,
    \quad
    \widetilde g^{(4)}(0)>0.
    \]
\end{enumerate}

Furthermore, degenerate minimizers occur exactly in the following four cases:

\begin{center}
\renewcommand{\arraystretch}{1.8}
\begin{tabular}{c|c|c|c}
\hline
Case & Parameter conditions & \(\mathcal S\) & \(\mathcal S_{\rm deg}\) \\
\hline
\((\mathrm{II}.1)\)
&
\(\begin{gathered}
r,s\in(1,2),\\
\displaystyle
\mu_1=\frac{2\kappa r}{p-2}
\left(\frac{2-s}{2-r}\right)^{\frac{s}{2}-1},
\\[0.4em]
\displaystyle
\mu_2=\frac{2\kappa s}{p-2}
\left(\frac{2-s}{2-r}\right)^{1-\frac{r}{2}}
\end{gathered}\)
&
\(\displaystyle
\left\{
\sqrt{\frac{2-s}{2-r}}
\right\}\)
&
\(\displaystyle
\left\{
\sqrt{\frac{2-s}{2-r}}
\right\}\)
\\
\hline
\((\mathrm{II}.2)\)
&
\(s=2,\, r>1,\, 0\le\mu_2<\mu_1=2\kappa,\,
  (r\ge2\ \text{or}\ \mu_2=0)\)
&
\(\{0\}\)
&
\(\{0\}\)
\\
\hline
\((\mathrm{II}.3)\)
&
\(r=2,\, s>1,\, 0\le\mu_1<\mu_2=2\kappa,\,
  (s\ge2\ \text{or}\ \mu_1=0)\)
&
\(\{\infty\}\)
&
\(\{\infty\}\)
\\
\hline
\((\mathrm{II}.4)\)
&
\(\min\{r,s\}=2<\max\{r,s\},\quad
\mu_1=\mu_2=2\kappa\)
&
\(\{0,\infty\}\)
&
\(\begin{cases}
\{0\}, & \text{if } s=2,\\
\{\infty\}, & \text{if } r=2.
\end{cases}\)
\\
\hline
\end{tabular}
\end{center}

In all other cases, every minimizer in \(\mathcal S\) is nondegenerate.
\end{lemma}

\begin{lemma}\label{lem2.7}\cite{ZZZ2026}
\begin{enumerate}
    \item
For any \(q\ge2\) and any \(m>0\), there exists \(C_1>0\) such that for any \(x,y\in\mathbb R\),
\[
|x+y|^q-|x|^q
\le q|x|^{q-2}xy
+\left(\frac{q(q-1)}2+m\right)|x|^{q-2}|y|^2
+C_1|y|^q .
\]
  \item
For any $\varsigma,\,\tau>1$ and any $m>0$, there exists a constant
$C=C(\varsigma,\tau,m)>0$ such that for any $x,z>0$ and
$y,w\in\mathbb R$, there holds
\[
\begin{aligned}
& |x+y|^\varsigma |z+w|^\tau-x^\varsigma z^\tau  \\
\le&
\varsigma x^{\varsigma-1}z^\tau y
+\tau x^\varsigma z^{\tau-1}w
+\big(\frac{\varsigma(\varsigma-1)}{2}+m\big)
x^{\varsigma-2}z^\tau y^2
+\big(\frac{\tau(\tau-1)}{2}+m\big)
x^\varsigma z^{\tau-2}w^2  \\
&+\varsigma\tau x^{\varsigma-1}z^{\tau-1}yw
+C\,\mathcal R_{\varsigma,\tau}(x,z;y,w),
\end{aligned}
\]
where $\chi_A$ denotes the characteristic function of the condition $A$, and
\[
\begin{aligned}
\mathcal R_{\varsigma,\tau}(x,z;y,w)
:=& |y|^\varsigma |w|^\tau +\chi_{\{\varsigma>2\}}\,z^\tau |y|^\varsigma
+\chi_{\{\tau>2\}}\,x^\varsigma |w|^\tau  \\
&+\chi_{\{\varsigma\ge 2\}}\,
x^{\varsigma-2}z^{\tau-1}|y|^2|w|
+\chi_{\{\tau\ge 2\}}\,
x^{\varsigma-1}z^{\tau-2}|y||w|^2  \\
&+\chi_{\{1<\varsigma<2\}}\,
x^{-1}z^\tau |y|^{\varsigma+1}
+\chi_{\{1<\tau<2\}}\,
x^\varsigma z^{-1}|w|^{\tau+1}  \\
&+\chi_{\{1<\varsigma<2\le \tau\}}\,
x^{-1}z^{\tau-1}|y|^{\varsigma+1}|w|  \\
&+\chi_{\{1<\tau<2\le \varsigma\}}\,
x^{\varsigma-1}z^{-1}|y||w|^{\tau+1}  \\
&+\chi_{\{1<\varsigma<2,\;1<\tau<2\}}
\big(
x^{\varsigma-1}z^{\frac{\tau-1}{2}}
|y||w|^{\frac{\tau+1}{2}}
+x^{\frac{\varsigma-1}{2}}z^{\tau-1}
|y|^{\frac{\varsigma+1}{2}}|w|
\big).
\end{aligned}
\]
\end{enumerate}
\end{lemma}

\section{The sharp system inequality}
In this section we prove Theorem \ref{th1.2}.
\begin{proof}[Proof of Theorem \ref{th1.2}]
For \((u,v)\in \mathcal X\), set
\[
a:=\left(\int_{\mathbb R^N}|x|^\ell |u|^p\,dx\right)^{1/p},
\qquad
b:=\left(\int_{\mathbb R^N}|x|^\ell |v|^p\,dx\right)^{1/p}.
\]
By the scalar sharp weighted fourth-order inequality, we have
\[
\|u\|_\alpha^2\ge S_0 a^2,\quad
\|v\|_\alpha^2\ge S_0 b^2.
\]
Hence
\begin{equation}\label{eq008}
\|(u,v)\|_{\mathcal X}^2
\ge S_0(a^2+b^2).
\end{equation}
On the other hand, by H\"older's inequality,
\[
\int_{\mathbb R^N}|x|^\ell |u|^r|v|^s\,dx
\le
\left(\int_{\mathbb R^N}|x|^\ell |u|^p\,dx\right)^{r/p}
\left(\int_{\mathbb R^N}|x|^\ell |v|^p\,dx\right)^{s/p}
=
a^r b^s.
\]
Consequently,
\[
P(u,v)
\le
\mu_1a^p+\mu_2b^p+p\kappa a^r b^s.
\]

We now compare the two homogeneous quantities. If \(a>0\), let \(t=b/a\). Then
\[
\frac{a^2+b^2}
{\big(\mu_1a^p+\mu_2b^p+p\kappa a^r b^s\big)^{2/p}}
=
\frac{1+t^2}
{\big(\mu_1+\mu_2t^p+p\kappa t^s\big)^{2/p}}
=g(t)
\ge m_*.
\]
If \(a=0\), then either \(b=0\), in which case the desired inequality is trivial, or
\[
\frac{a^2+b^2}
{\big(\mu_1a^p+\mu_2b^p+p\kappa a^r b^s\big)^{2/p}}
=
\mu_2^{-2/p}
=
g(\infty)
\ge m_*.
\]
Therefore, in all cases,
\[
a^2+b^2
\ge
m_*
\big(\mu_1a^p+\mu_2b^p+p\kappa a^r b^s\big)^{2/p}
\ge
m_* P(u,v)^{2/p}.
\]
Combining \eqref{eq008} gives
\[
\|(u,v)\|_{\mathcal X}^2
\ge
m_*S_0 P(u,v)^{2/p}.
\]
Thus the system inequality holds with the constant
\[
S_{r,s,\mu_1,\mu_2,\kappa}=m_*S_0.
\]

It remains to prove that this constant is sharp and to characterize the equality cases.
Let \(t_0\in \mathcal S\cap[0,\infty)\), let \(\lambda>0\), and let \(A\in\mathbb R\). Taking
\[
(u,v)=A(W_\lambda,\sigma t_0W_\lambda),
\qquad \sigma\in\{-1,1\},
\]
we have equality in the scalar inequalities and in H\"older's inequality. Moreover,
\[
\frac{\|(u,v)\|_{\mathcal X}^2}{P(u,v)^{2/p}}
=
S_0 g(t_0)
=
m_*S_0.
\]
If \(\infty\in\mathcal S\), the same conclusion follows by taking
\[
(u,v)=(0,AW_\lambda).
\]
Therefore the constant \(m_*S_0\) is optimal.

Conversely, suppose that equality holds in the system inequality for some
\((u,v)\not=(0,0)\). Then equality must hold in each of the inequalities used above. In
particular, the ratio \(t=b/a\), with the convention \(t=\infty\) if \(a=0\), must belong to
\(\mathcal S\).

Assume first that \(t\in(0,\infty)\). Then \(a,b>0\). Equality in the scalar inequalities implies
that there exist \(c,d\in\mathbb R\setminus\{0\}\) and \(\lambda,\mu>0\) such that
\[
u=cW_\lambda,\qquad v=dW_\mu.
\]
Equality in H\"older's inequality gives
\[
\frac{|u|^p}{a^p}=\frac{|v|^p}{b^p}
\quad\text{a.e. in }\mathbb R^N.
\]
Hence \(W_\lambda\) and \(W_\mu\) are proportional. Since the scalar optimizers are normalized
up to the same dilation family, this forces \(\lambda=\mu\). Moreover,
\[
\frac{|d|}{|c|}=\frac{b}{a}=t.
\]
Thus, for some \(\sigma\in\{-1,1\}\),
\[
(u,v)=c(W_\lambda,\sigma tW_\lambda).
\]

If \(t=0\), then \(b=0\), and therefore \(v=0\). Equality in the scalar inequality for \(u\)
gives
\[
(u,v)=(cW_\lambda,0).
\]
Similarly, if \(t=\infty\),
\[
(u,v)=(0,cW_\lambda).
\]

Consequently, the full set of equality cases is
\[
\left\{
A(W_\lambda,\sigma tW_\lambda):
A\in\mathbb R\setminus\{0\},\ \lambda>0,\ t\in \mathcal S\cap[0,\infty),\sigma\in\{-1,1\}
\right\},
\]
together with the branch
\[
\{(0,AW_\lambda):A\in\mathbb R\setminus\{0\},\lambda>0\}
\]
if \(\infty\in\mathcal S\). This is precisely \(\mathcal M\).
In particular, among nonnegative pairs, one takes \(A>0\) and \(\sigma=1\) whenever
\(t\in[0,\infty)\).

If \((r,s,\mu_1,\mu_2)=(2,2,2\kappa,2\kappa)\), then
\[
P(u,v)=
2\kappa
\int_{\mathbb R^N}|x|^\ell
\bigl(|u|^2+|v|^2\bigr)^2\,dx,\quad g(t)
=\frac{1}{\sqrt{2\kappa}},
\]
and hence
\[
\|(u,v)\|_{\mathcal X}^2
\geq
S_0
\left(
\int_{\mathbb R^N}|x|^\ell
\bigl(|u|^2+|v|^2\bigr)^2\,dx
\right)^{1/2}.
\]
Equality holds for a nonzero pair \((u,v)\) if and only if
\[
(u,v)=(aW_\lambda,bW_\lambda)
\]
for some
\(
(a,b)\in\mathbb R^2\setminus\{(0,0)\}\)
and
\(\lambda>0.
\)
\end{proof}

\section{Qualitative stability and orthogonal projection}
\label{sec:qualitative}
In this section, we establish the compactness of normalized minimizing
sequences and a qualitative stability result near the minimizing
manifold. These results allow us to reduce the proof of the sharp
stability estimate to a sufficiently small neighborhood of the
minimizing manifold.

Throughout this section, we assume that one of \eqref{case-A} and \eqref{case-B} holds and that
\[
(r,s,\mu_1,\mu_2)\neq(2,2,2\kappa,2\kappa).
\]
Thus the set $\mathcal S$ of minimizers of $g$ is finite.

Since $\operatorname{dist}\bigl((u,v),\mathcal M\bigr)^{2}\leq
\|(u,v)\|_{\mathcal X}^2$, if
\(
\delta(u,v)
\geq
\varepsilon_{0}
\bigl(
\|(u,v)\|_{\mathcal X}^2
\bigr),
\)
then
\[
\left(
\frac{\delta(u,v)}
{\|(u,v)\|_{\mathcal X}^2}
\right)^{\iota}
\geq
\varepsilon_{0}^{\iota}
\frac{
\operatorname{dist}\bigl((u,v),\mathcal M\bigr)^{2}
}{
\|(u,v)\|_{\mathcal X}^2
}.
\]
Hence, after changing the stability constant if necessary, the
large-deficit case is immediate. It is therefore enough to study the
small-deficit regime.

For $t\in \mathcal S\cap[0,\infty)$ and
$\sigma\in\{-1,1\}$, define
\[
\eta_{t,\sigma}(u,v)
:=
\inf_{A\in\mathbb R\setminus\{0\},\ \lambda>0}
\left\{
\|u-AW_{\lambda}\|_{\alpha}^{2}
+
\|v-\sigma tAW_{\lambda}\|_{\alpha}^{2}
\right\}.
\]
For $t=\infty$, define
\[
\eta_{\infty}(u,v)
:=
\inf_{A\in\mathbb R\setminus\{0\},\ \lambda>0}
\left\{
\|u\|_{\alpha}^{2}
+
\|v-AW_{\lambda}\|_{\alpha}^{2}
\right\}.
\]

We first prove that the above projection problems admit minimizers
whenever the distance is strictly smaller than the distance to the
origin.

\begin{lemma}
\label{lem:projection-attainment}
Let $t\in \mathcal S\cap[0,\infty)$ and
$\sigma\in\{-1,1\}$. Suppose that
$(u,v)\in\mathcal X$ satisfies
\[
\eta_{t,\sigma}(u,v)
<
\|(u,v)\|_{\mathcal X}^2.
\]
Then there exist
\[
A\in\mathbb R\setminus\{0\},
\quad
\lambda\in(0,\infty),
\]
such that
\[
\eta_{t,\sigma}(u,v)
=
\|u-AW_{\lambda}\|_{\alpha}^{2}
+
\|v-\sigma tAW_{\lambda}\|_{\alpha}^{2}.
\]

Similarly, if
\[
\eta_{\infty}(u,v)
<
\|(u,v)\|_{\mathcal X}^2,
\]
then there exist
\[
A\in\mathbb R\setminus\{0\},
\quad
\lambda\in(0,\infty),
\]
such that
\[
\eta_{\infty}(u,v)
=
\|u\|_{\alpha}^{2}
+
\|v-AW_{\lambda}\|_{\alpha}^{2}.
\]
\end{lemma}

\begin{proof}
We first consider $t\in \mathcal S\cap[0,\infty)$. Set
\(
E:=\|(u,v)\|_{\mathcal X}^2
\)
and
\[
L(\lambda)
:=
\langle u,W_{\lambda}\rangle_{\alpha}
+
\sigma t
\langle v,W_{\lambda}\rangle_{\alpha},
\]
where
\[
\langle f,h\rangle_{\alpha}
:=
\int_{\mathbb R^{N}}
|x|^{\alpha}\Delta f\,\Delta h\,dx.
\]
For fixed $\lambda>0$, we have
\begin{align*}
\|u-AW_{\lambda}\|_{\alpha}^{2}
+
\|v-\sigma tAW_{\lambda}\|_{\alpha}^{2}
=E-2AL(\lambda)+
A^{2}(1+t^{2})S_0.
\end{align*}
The unique minimizer with respect to $A$ is
\[
A_{\lambda}
=
\frac{L(\lambda)}
{(1+t^{2})S_0}.
\]
The corresponding minimum is
\[
E-
\frac{|L(\lambda)|^{2}}
{(1+t^{2})S_0}.
\]
Consequently,
\[
\eta_{t,\sigma}(u,v)
=
E-
\frac{1}
{(1+t^{2})S_0}
\sup_{\lambda>0}|L(\lambda)|^{2}.
\]
The strict inequality
\[
\eta_{t,\sigma}(u,v)<E
\]
therefore implies
\[
\sup_{\lambda>0}|L(\lambda)|>0.
\]

By Lemma \ref{lem:weak-bubble}, there exist $\lambda_0\in(0,\infty)$, such that
\[
|L(\lambda_0)|=\sup_{\lambda>0}|L(\lambda)|
\]
and
\[
A_{0}
=
\frac{L(\lambda_{0})}
{(1+t^{2})S_0}.
\]
Then $A_{0}\neq0$ and
\[
\eta_{t,\sigma}(u,v)
=
\|u-A_{0}W_{\lambda_{0}}\|_{\alpha}^{2}
+
\|v-\sigma tA_{0}W_{\lambda_{0}}\|_{\alpha}^{2}.
\]

For $t=\infty$, define
\[
L_{\infty}(\lambda)
:=
\langle v,W_{\lambda}\rangle_{\alpha}.
\]
For fixed $\lambda$, the unique minimizer in $A$ is
\[
A_{\lambda}
=
\frac{L_{\infty}(\lambda)}
{S_0}.
\]
The same argument shows that
$|L_{\infty}|$ attains its positive supremum at some
$\lambda_{0}\in(0,\infty)$, which proves the assertion.
\end{proof}

\begin{lemma}\label{lem:overlap}
Let
\[
I(\gamma):=\int_{\mathbb R^N}|x|^\ell W_1^rW_\gamma^s\,dx,\quad \gamma>0.
\]
Then \(I(\gamma)\le1\), equality holds if and only if \(\gamma=1\), and
\[
I(\gamma)\longrightarrow0
\quad\text{as }\gamma\to0\text{ or }\gamma\to\infty.
\]
Consequently, if \(I(\gamma_n)\to1\), then \(\gamma_n\to1\).
\end{lemma}

\begin{proof}
By H\"older's inequality and the normalization
\[
\int_{\mathbb R^N}|x|^\ell W_\gamma^p\,dx=1,
\]
we have
\[
I(\gamma)\le
\left(\int_{\mathbb R^N}|x|^\ell W_1^p\,dx\right)^{r/p}
\left(\int_{\mathbb R^N}|x|^\ell W_\gamma^p\,dx\right)^{s/p}
=1.
\]
The equality condition in H\"older's inequality implies that \(W_1^p\) and
\(W_\gamma^p\) are linearly dependent. Since they have the same normalization,
this gives \(W_\gamma=W_1\), and hence \(\gamma=1\).

It remains to show the vanishing at the two ends. Let \(\gamma_n\to0\) or
\(\gamma_n\to\infty\). By the dilation structure of the bubbles, the family
\(W_{\gamma_n}^{s}\) converges weakly to zero in
\(L^{p/s}(\mathbb R^N,|x|^\ell dx)\). Since
\(W_1^r\in L^{p/r}(\mathbb R^N,|x|^\ell dx)\), we obtain
\(I(\gamma_n)\to0\). This weak convergence can also be checked directly by
splitting the integral on a fixed ball and its complement and using the
explicit decay of \(W_\gamma\) in Cases A and B. Since \(I\) is continuous on
compact subintervals of \((0,\infty)\), the final assertion follows.
\end{proof}

\begin{proposition}
\label{prop:system-compactness}
Let
\(
(u_{n},v_{n})\subset\mathcal X
\)
satisfy
\(
\|(u_n,v_n)\|_{\mathcal X}^2
=1
\)
and
\(
\delta(u_{n},v_{n})\longrightarrow0.
\)
Then, after passing to a subsequence, there exists
$t_{0}\in \mathcal S$ such that one of the following alternatives holds.

\begin{enumerate}
\item If $t_{0}\in[0,\infty)$, then there exist
$\sigma\in\{-1,1\}$, $A_{n}\in\mathbb R$, and
$\lambda_{n}>0$ such that
\[
\|u_{n}-A_{n}W_{\lambda_{n}}\|_{\alpha}^{2}
+
\|v_{n}-\sigma t_{0}A_{n}W_{\lambda_{n}}\|_{\alpha}^{2}
\longrightarrow0.
\]

\item If $t_{0}=\infty$, then there exist
$A_{n}\in\mathbb R$ and $\lambda_{n}>0$ such that
\[
\|u_{n}\|_{\alpha}^{2}
+
\|v_{n}-A_{n}W_{\lambda_{n}}\|_{\alpha}^{2}
\longrightarrow0.
\]
\end{enumerate}

In particular,
\[
\operatorname{dist}
\bigl((u_{n},v_{n}),\mathcal M\bigr)
\longrightarrow0.
\]
\end{proposition}

\begin{proof}
For simplicity, write
\[
P_{n}:=P(u_{n},v_{n}).
\]
Since
\[
\delta(u_{n},v_{n})
=
1-m^{*}S_{0}P_{n}^{2/p}
\longrightarrow0,
\]
we obtain
\[
P_{n}^{2/p}
\longrightarrow
\frac{1}{m^{*}S_{0}}.
\]
In particular, $P_{n}$ is bounded away from zero and infinity.

Define
\[
a_{n}
:=
\left(
\int_{\mathbb R^{N}}
|x|^{\ell}|u_{n}|^{p}\,dx
\right)^{1/p},
\qquad
b_{n}
:=
\left(
\int_{\mathbb R^{N}}
|x|^{\ell}|v_{n}|^{p}\,dx
\right)^{1/p},
\]
and
\[
F_{n}
:=
\mu_{1}a_{n}^{p}
+
\mu_{2}b_{n}^{p}
+
p\kappa a_{n}^{r}b_{n}^{s}.
\]
By H\"older's inequality,
\[
P_{n}
\leq
F_{n}.
\]
On the other hand, the scalar weighted Caffarelli--Kohn--Nirenberg
inequality gives
\[
1
\geq
S_{0}(a_{n}^{2}+b_{n}^{2}).
\]

Let
\[
t_{n}:=
\frac{b_{n}}{a_{n}}
\]
when $a_{n}>0$, with the convention $t_{n}=\infty$ when
$a_{n}=0$. We thus have
\begin{equation}
\frac{1}{P_{n}^{2/p}}
\geq
\frac{1}{F_{n}^{2/p}}
\geq
S_{0}
\frac{a_{n}^{2}+b_{n}^{2}}{F_{n}^{2/p}}
=
S_{0}g(t_{n})
\geq
S_{0}m^{*}.
\label{eq:chain-compactness}
\end{equation}
The left-hand side converges to $S_{0}m^{*}$. Hence every
nonnegative gap in \eqref{eq:chain-compactness} tends to zero. In
particular,
\[
g(t_{n})\longrightarrow m^{*}.
\]
Since $\mathcal S$ is finite, after passing to a subsequence,
\[
t_{n}\longrightarrow t_{0}
\quad\text{for some }t_{0}\in \mathcal S,
\]
where convergence is understood in $[0,\infty]$.

Moreover,
\[
\frac{1}{F_{n}^{2/p}}
\longrightarrow S_{0}m^{*}.
\]
Hence $F_{n}$ is bounded away from zero and infinity. From
\[
\frac{1-S_{0}(a_{n}^{2}+b_{n}^{2})}{F_{n}^{2/p}}
=
\frac{1}{F_{n}^{2/p}}
-
S_{0}
\frac{a_{n}^{2}+b_{n}^{2}}{F_{n}^{2/p}}
\longrightarrow0,
\]
we obtain
\[
1-S_{0}(a_{n}^{2}+b_{n}^{2})
\longrightarrow0.
\]
Therefore
\begin{equation}
\|u_{n}\|_{\alpha}^{2}-S_{0}a_{n}^{2}
\longrightarrow0,
\quad
\|v_{n}\|_{\alpha}^{2}-S_{0}b_{n}^{2}
\longrightarrow0.
\label{eq:scalar-deficits-vanish}
\end{equation}

The first gap in \eqref{eq:chain-compactness} gives
\[
\frac{1}{P_{n}^{2/p}}
-
\frac{1}{F_{n}^{2/p}}
\longrightarrow0.
\]
Since $P_{n}$ and $F_{n}$ stay in a compact subset of $(0,\infty)$,
it follows that
\[
F_{n}-P_{n}\longrightarrow0.
\]
Hence
\begin{equation}
a_{n}^{r}b_{n}^{s}
-
\int_{\mathbb R^{N}}
|x|^{\ell}|u_{n}|^{r}|v_{n}|^{s}\,dx
\longrightarrow0.
\label{eq:holder-deficit-vanish}
\end{equation}

We now distinguish three cases.

\medskip
\noindent
\textbf{Case 1: $t_{0}\in(0,\infty)$.}

Since $t_{n}\to t_{0}$ and $F_{n}$ is bounded away from zero,
both $a_{n}$ and $b_{n}$ are bounded away from zero. By
\eqref{eq:scalar-deficits-vanish} and the scalar stability estimate
of Theorem~\ref{thm:scalar-stability}, there exist
\[
c_{n},d_{n}\in\mathbb R\setminus\{0\},
\quad
\lambda_{n},\mu_{n}>0,
\]
such that
\begin{equation}\label{eq:two-scalar-bubbles}
u_{n}-c_{n}W_{\lambda_{n}}
\longrightarrow0,
\quad
v_{n}-d_{n}W_{\mu_{n}}
\longrightarrow0
\quad\text{in }X_{\mathrm{rad}}.
\end{equation}
By the continuous embedding
\(
X_{\mathrm{rad}}
\hookrightarrow
L^{p}(\mathbb R^{N},|x|^{\ell}dx)
\)
and the normalization \eqref{scalar-normalization}, we have
\begin{equation}
a_{n}=|c_{n}|+o(1),
\quad
b_{n}=|d_{n}|+o(1).
\label{eq:amplitudes-lp}
\end{equation}

The convergence in \eqref{eq:two-scalar-bubbles}, together with
H\"older's inequality, yields
\begin{equation}\label{eq0001}
\int_{\mathbb R^{N}}
|x|^{\ell}|u_{n}|^{r}|v_{n}|^{s}\,dx
=|c_{n}|^{r}|d_{n}|^{s}
\int_{\mathbb R^{N}}
|x|^{\ell}
W_{\lambda_{n}}^{r}W_{\mu_{n}}^{s}\,dx
+o(1).
\end{equation}
Combining this with
\eqref{eq:holder-deficit-vanish}, \eqref{eq:amplitudes-lp} and \eqref{eq0001}, and using that
$|c_{n}|$ and $|d_{n}|$ are bounded away from zero, we obtain
\begin{equation}
\int_{\mathbb R^{N}}
|x|^{\ell}
W_{\lambda_{n}}^{r}W_{\mu_{n}}^{s}\,dx
\longrightarrow1.
\label{eq:overlap-one}
\end{equation}

Set
\[
\gamma_{n}:=\frac{\mu_{n}}{\lambda_{n}}.
\]
By the critical scaling of $W_{\lambda}$,
\[
\int_{\mathbb R^{N}}
|x|^{\ell}
W_{\lambda_{n}}^{r}W_{\mu_{n}}^{s}\,dx
=
I(\gamma_{n}),
\]
Hence \eqref{eq:overlap-one} and Lemma~\ref{lem:overlap} imply
\[
\gamma_{n}\longrightarrow1.
\]
Consequently,
\[
\|W_{\mu_{n}}-W_{\lambda_{n}}\|_{\alpha}
\longrightarrow0,
\]
and hence
\[
v_{n}-d_{n}W_{\lambda_{n}}
\longrightarrow0
\quad\text{in }X_{\mathrm{rad}}.
\]

Since
\[
\frac{b_{n}}{a_{n}}\longrightarrow t_{0},
\]
relation \eqref{eq:amplitudes-lp} gives
\[
\frac{|d_{n}|}{|c_{n}|}\longrightarrow t_{0}.
\]
After passing to a further subsequence, the sign of $d_{n}/c_{n}$
is constant. Then
\[
d_{n}-\sigma t_{0}c_{n}\longrightarrow0
\]
for some \(\sigma\in\{-1,1\}\). Therefore
\[
\|u_{n}-c_{n}W_{\lambda_{n}}\|_{\alpha}^{2}
+
\|v_{n}-\sigma t_{0}c_{n}W_{\lambda_{n}}\|_{\alpha}^{2}
\longrightarrow0.
\]

\medskip
\noindent
\textbf{Case 2: $t_{0}=0$.}

In this case,
\[
\frac{b_{n}}{a_{n}}\longrightarrow0.
\]
Since $F_{n}$ is bounded away from zero, $a_{n}$ is bounded away
from zero, whereas
\[
b_{n}\longrightarrow0.
\]
By \eqref{eq:scalar-deficits-vanish},
\[
\|v_{n}\|_{\alpha}^{2}
=
S_{0}b_{n}^{2}+o(1)
\longrightarrow0.
\]
The first relation in \eqref{eq:scalar-deficits-vanish}, together
with Theorem~\ref{thm:scalar-stability}, gives $c_{n}\in\mathbb R\setminus\{0\}$ and
$\lambda_{n}>0$ such that
\[
u_{n}-c_{n}W_{\lambda_{n}}
\longrightarrow0
\quad\text{in }X_{\mathrm{rad}}.
\]
It follows that
\[
\|u_{n}-c_{n}W_{\lambda_{n}}\|_{\alpha}^{2}
+
\|v_{n}\|_{\alpha}^{2}
\longrightarrow0.
\]
This is precisely the convergence to the semi-trivial branch
corresponding to $t_{0}=0$.

\medskip
\noindent
\textbf{Case 3: $t_{0}=\infty$.}

Here
\[
\frac{b_{n}}{a_{n}}\longrightarrow\infty.
\]
Since $F_{n}$ is bounded away from zero and
$a_{n},b_{n}$ are bounded by the scalar inequality, we have
\[
a_{n}\longrightarrow0.
\]
By
\eqref{eq:scalar-deficits-vanish},
\[
\|u_{n}\|_{\alpha}^{2}
=
S_{0}a_{n}^{2}+o(1)
\longrightarrow0.
\]
Applying Lemma~\ref{lem:scalar-cylinder-compactness} to the second component, we find
$d_{n}\in\mathbb R\setminus\{0\}$ and $\mu_{n}>0$ such that
\[
v_{n}-d_{n}W_{\mu_{n}}
\longrightarrow0
\quad\text{in }X_{\mathrm{rad}}.
\]
Thus
\[
\|u_{n}\|_{\alpha}^{2}
+
\|v_{n}-d_{n}W_{\mu_{n}}\|_{\alpha}^{2}
\longrightarrow0,
\]
which is the desired convergence to the branch $t_{0}=\infty$.
\end{proof}

We now convert the sequential compactness result into a uniform
qualitative stability statement.

\begin{theorem}
\label{thm:qualitative-stability}
For every $\varepsilon_{1}>0$, there exists
$\varepsilon_{0}>0$ such that the following holds. If
$(u,v)\in\mathcal X$ satisfies
\[
0<\delta(u,v)
<
\varepsilon_{0}
\bigl(
\|u\|_{\alpha}^{2}+\|v\|_{\alpha}^{2}
\bigr),
\]
then the following assertion holds. There exists $t_{0}\in\mathcal S$. If $t_{0}<\infty$, then for some $\sigma\in\{-1,1\}$ there is a minimizing pair $(A,\lambda)\in\mathbb R\times(0,\infty)$ for $\eta_{t_{0},\sigma}$ such that $A\neq0$ and
\[
\eta_{t_{0},\sigma}(u,v)
=
\|u-AW_{\lambda}\|_{\alpha}^{2}
+
\|v-\sigma t_{0}AW_{\lambda}\|_{\alpha}^{2}
<
\varepsilon_{1}A^{2}.
\]
If $t_{0}=\infty$, then there is a minimizing pair $(A,\lambda)\in\mathbb R\times(0,\infty)$ for $\eta_\infty$ such that $A\neq0$ and
\[
\eta_{\infty}(u,v)
=
\|u\|_{\alpha}^{2}
+
\|v-AW_{\lambda}\|_{\alpha}^{2}
<
\varepsilon_{1}A^{2}.
\]

In particular, for every sequence
$(u_{n},v_{n})\subset\mathcal X\setminus\{(0,0)\}$,
\[
\frac{\delta(u_{n},v_{n})}
{\|(u_n,v_n)\|_{\mathcal X}^2}
\longrightarrow0
\]
implies
\[
\frac{
\operatorname{dist}
\bigl((u_{n},v_{n}),\mathcal M\bigr)^{2}
}{
\|(u_n,v_n)\|_{\mathcal X}^2
}
\longrightarrow0.
\]
\end{theorem}

\begin{proof}
Suppose, by contradiction, that the assertion is false. Then there
exist $\varepsilon_{1}>0$ and a sequence
\[
(u_{n},v_{n})\subset\mathcal X
\]
such that, by the homogeneity of the deficit and the projection
problems,
\[
\|(u,v)\|_{\mathcal X}^2=1,
\quad
0<\delta(u_{n},v_{n})\leq\frac1n,
\]
but the conclusion fails for every minimizing branch. More precisely,
for every $t\in \mathcal S\cap[0,\infty)$, every
$\sigma\in\{-1,1\}$, and every pair $(A,\lambda)$ attaining
$\eta_{t,\sigma}(u_{n},v_{n})$, one has
\[
\eta_{t,\sigma}(u_{n},v_{n})
\geq
\varepsilon_{1}A^{2}.
\]
Likewise, for every pair $(A,\lambda)$ attaining
$\eta_{\infty}(u_{n},v_{n})$, one has
\[
\eta_{\infty}(u_{n},v_{n})
\geq
\varepsilon_{1}A^{2}.
\]

By Proposition~\ref{prop:system-compactness}, after passing to a
subsequence, there exists $t_{0}\in \mathcal S$ for which the sequence
approaches the corresponding branch.

Assume first that $t_{0}\in[0,\infty)$. Then there exist
$\sigma\in\{-1,1\}$, $a_{n}\in\mathbb R$, and
$\rho_{n}>0$ such that
\[
\|u_{n}-a_{n}W_{\rho_{n}}\|_{\alpha}^{2}
+
\|v_{n}-\sigma t_{0}a_{n}W_{\rho_{n}}\|_{\alpha}^{2}
\longrightarrow0.
\]
Consequently,
\[
0
\leq
\eta_{t_{0},\sigma}(u_{n},v_{n})
\leq
\|u_{n}-a_{n}W_{\rho_{n}}\|_{\alpha}^{2}
+
\|v_{n}-\sigma t_{0}a_{n}W_{\rho_{n}}\|_{\alpha}^{2}
\longrightarrow0.
\]
Set
\[
\eta_{n}
:=
\eta_{t_{0},\sigma}(u_{n},v_{n}).
\]
Lemma~\ref{lem:projection-attainment}
therefore implies that $\eta_{n}$ is attained by some
\[
(\widehat A_{n},\widehat\lambda_{n})
\in
\bigl(\mathbb R\setminus\{0\}\bigr)\times(0,\infty).
\]
Thus
\[
\eta_{n}
=
\|u_{n}-\widehat A_{n}W_{\widehat\lambda_{n}}\|_{\alpha}^{2}
+
\|v_{n}-\sigma t_{0}\widehat A_{n}
W_{\widehat\lambda_{n}}\|_{\alpha}^{2}.
\]
By the reverse triangle inequality,
\[
|\widehat A_{n}|\sqrt{(1+t_{0}^{2})S_0}
=\|\widehat A_{n}
\bigl(
W_{\widehat\lambda_{n}},
\sigma t_{0}W_{\widehat\lambda_{n}}
\bigr)\|_{\mathcal X}
\geq
1-\sqrt{\eta_{n}}.
\]
Therefore
\[
\widehat A_{n}^{\,2}
\geq
\frac{(1-\sqrt{\eta_{n}})^{2}}
{(1+t_{0}^{2})S_0}.
\]
In particular, there exists $c_{0}>0$ such that
\[
\widehat A_{n}^{\,2}\geq c_{0}
\]
for all sufficiently large $n$. It follows that
\[
0
\leq
\frac{\eta_{n}}{\widehat A_{n}^{\,2}}
\leq
\frac{\eta_{n}}{c_{0}}
\longrightarrow0.
\]
Hence, for sufficiently large $n$,
\[
\eta_{t_{0},\sigma}(u_{n},v_{n})
=
\eta_{n}
<
\varepsilon_{1}\widehat A_{n}^{\,2},
\]
contradicting the assumed failure of the conclusion.

If $t_{0}=\infty$, Proposition~\ref{prop:system-compactness}
gives $a_{n}\in\mathbb R$ and $\rho_{n}>0$ such that
\[
\|u_{n}\|_{\alpha}^{2}
+
\|v_{n}-a_{n}W_{\rho_{n}}\|_{\alpha}^{2}
\longrightarrow0.
\]
Thus
\[
\eta_{\infty}(u_{n},v_{n})
\longrightarrow0.
\]
For sufficiently large $n$, Lemma~\ref{lem:projection-attainment}
gives a minimizing pair
$(\widehat A_{n},\widehat\lambda_{n})$ satisfying
\[
\eta_{\infty}(u_{n},v_{n})
=
\|u_{n}\|_{\alpha}^{2}
+
\|v_{n}-\widehat A_{n}W_{\widehat\lambda_{n}}\|_{\alpha}^{2}.
\]
Similarly,
\[
|\widehat A_{n}|\sqrt{S_0}
=
\|(0,\widehat A_{n}W_{\widehat\lambda_{n}})\|_{\mathcal X}
\geq
1-\sqrt{\eta_{\infty}(u_{n},v_{n})}.
\]
Hence $\widehat A_{n}^{\,2}$ is bounded away from zero, and
\[
\frac{\eta_{\infty}(u_{n},v_{n})}
{\widehat A_{n}^{\,2}}
\longrightarrow0.
\]
This again contradicts the assumed failure of the theorem.

It remains to prove the final sequential statement. Let
\[
E_{n}:=\|(u_n,v_n)\|_{\mathcal X}^2>0
\]
and define
\[
\widetilde u_{n}
:=
\frac{u_{n}}{\sqrt{E_{n}}},
\quad
\widetilde v_{n}
:=
\frac{v_{n}}{\sqrt{E_{n}}}.
\]
By homogeneity,
\[\|(\widetilde u_n,\widetilde v_n)\|_{\mathcal X}^2=1\]
and
\[
\delta(\widetilde u_{n},\widetilde v_{n})
=
\frac{\delta(u_{n},v_{n})}{E_{n}}
\longrightarrow0.
\]
Proposition~\ref{prop:system-compactness} therefore gives
\[
\operatorname{dist}
\bigl(
(\widetilde u_{n},\widetilde v_{n}),
\mathcal M
\bigr)
\longrightarrow0.
\]
Since $\mathcal M$ is a cone,
\[
\operatorname{dist}
\bigl(
(\widetilde u_{n},\widetilde v_{n}),
\mathcal M
\bigr)
=
\frac{
\operatorname{dist}
\bigl((u_{n},v_{n}),\mathcal M\bigr)
}{
\sqrt{E_{n}}
}.
\]
The desired conclusion follows.
\end{proof}

The minimizing pairs obtained above provide the orthogonality
conditions required in the local analysis of the next section.

\begin{corollary}\label{cor:projection-orthogonality}
Let $t_{0}\in \mathcal S\cap[0,\infty)$ and
$\sigma\in\{-1,1\}$, and suppose that
$(A,\lambda)$ is a minimizing pair for
$\eta_{t_{0},\sigma}(u,v)$ with $A\neq0$. Set
\[
\phi:=u-AW_{\lambda},
\quad
\psi:=v-\sigma t_{0}AW_{\lambda}.
\]
Then
\[
\langle\phi,W_{\lambda}\rangle_{\alpha}
+
\sigma t_{0}
\langle\psi,W_{\lambda}\rangle_{\alpha}
=0
\]
and
\[
\langle\phi,\partial_{\lambda}W_{\lambda}\rangle_{\alpha}
+
\sigma t_{0}
\langle\psi,\partial_{\lambda}W_{\lambda}\rangle_{\alpha}
=0.
\]

If $t_{0}=\infty$ and $(A,\lambda)$ minimizes
$\eta_{\infty}(u,v)$, then
\[
\langle v-AW_{\lambda},W_{\lambda}\rangle_{\alpha}
=0
\]
and
\[
\langle v-AW_{\lambda},
\partial_{\lambda}W_{\lambda}\rangle_{\alpha}
=0.
\]
\end{corollary}

\begin{proof}
For $t_{0}\in[0,\infty)$, consider
\[
\mathcal F(A,\lambda)
:=
\|u-AW_{\lambda}\|_{\alpha}^{2}
+
\|v-\sigma t_{0}AW_{\lambda}\|_{\alpha}^{2}.
\]
Since $(A,\lambda)$ is an interior minimizing pair and $A\neq0$,
differentiating with respect to $A$ gives
\[
\langle\phi,W_{\lambda}\rangle_{\alpha}
+
\sigma t_{0}
\langle\psi,W_{\lambda}\rangle_{\alpha}
=0.
\]
Differentiating with respect to $\lambda$ gives
\[
A\left(
\langle\phi,\partial_{\lambda}W_{\lambda}\rangle_{\alpha}
+
\sigma t_{0}
\langle\psi,\partial_{\lambda}W_{\lambda}\rangle_{\alpha}
\right)
=0.
\]
Since $A\neq0$, the second orthogonality condition follows.

The case $t_{0}=\infty$ is obtained in exactly the same way by
differentiating
\[
\|u\|_{\alpha}^{2}
+
\|v-AW_{\lambda}\|_{\alpha}^{2}
\]
with respect to $A$ and $\lambda$.
\end{proof}

We finally record a local geometric estimate which will be used in the local
analysis of Section~\ref{sec:local-analysis}. It identifies the distance to a fixed synchronized
branch with the size of the normal components in suitable local coordinates. For \(t_0\in\mathcal S\cap(0,\infty)\), set
\[
e_0:=\frac{(1,t_0)}{\sqrt{1+t_0^2}},
\quad
e_1:=\frac{(-t_0,1)}{\sqrt{1+t_0^2}}.
\]
Then \(e_0\) is the synchronized direction and \(e_1\) is the transverse
ratio-changing direction.

\begin{lemma}\label{lem:local-distance}
Let \(t_0\in\mathcal S\cap(0,\infty)\) and let
\[
\mathcal M_{t_0}:=\{A(W_\lambda,t_0W_\lambda):A\in\mathbb R,\lambda>0\}.
\]
Assume that
\[
(u,v)=A(W_\lambda,t_0W_\lambda)+e_0\xi+e_1(\rho W_\lambda+\zeta),
\]
where \(A>0\), \(\lambda>0\), and
\[
\langle\xi,W_\lambda\rangle_\alpha=
\langle\xi,\partial_\lambda W_\lambda\rangle_\alpha=
\langle\zeta,W_\lambda\rangle_\alpha=0.
\]
If
\[
\|\xi\|_\alpha+\|\zeta\|_\alpha+|\rho|\le \varepsilon A
\]
with \(\varepsilon>0\) sufficiently small, then
\[
\operatorname{dist}((u,v),\mathcal M_{t_0})^2
\asymp
\|\xi\|_\alpha^2+\|\zeta\|_\alpha^2+\rho^2.
\]
The constants depend only on the fixed parameters and on \(t_0\). The analogous statement holds for the endpoint branches after replacing \(\mathcal M_{t_0}\) by the corresponding semi-trivial branch.
\end{lemma}

\begin{proof}
The upper bound follows by taking the fixed point
\(A(W_\lambda,t_0W_\lambda)\) as a competitor. We prove the lower bound. The
true tangent space to \(\mathcal M_{t_0}\) at \(A(W_\lambda,t_0W_\lambda)\) is
\[
        T_{A,\lambda}\mathcal M_{t_0}
        =\spanop\{e_0W_\lambda,e_0\partial_\lambda W_\lambda\}.
\]
The direction \(e_1W_\lambda\) changes the ratio of the two components and is
not tangent to \(\mathcal M_{t_0}\). Thus it is a normal ratio direction and
must be counted in the distance to the fixed branch. Since \(e_0\perp e_1\) in
\(\mathbb R^2\), the imposed orthogonality conditions give
\[
        e_0\xi+e_1(\rho W_\lambda+\zeta)
        \perp T_{A,\lambda}\mathcal M_{t_0},
\]
and
\[
        \|e_1(\rho W_\lambda+\zeta)\|_{\mathcal X}^2
        =\rho^2\|W_\lambda\|_\alpha^2+\|\zeta\|_\alpha^2.
\]
The map \((A,\lambda)\mapsto A(W_\lambda,t_0W_\lambda)\) is \(C^1\). Moreover,
the Gram matrix of the tangent vectors \(e_0W_\lambda\) and
\(e_0\partial_\lambda W_\lambda\) is nonsingular after the usual normalization of
\(\partial_\lambda W_\lambda\), because \(\partial_\lambda W_\lambda\not\equiv0\) and
\(\langle W_\lambda,\partial_\lambda W_\lambda\rangle_\alpha=0\), the latter
identity following from differentiating \(\|W_\lambda\|_\alpha^2=S_0\). Hence,
in a sufficiently small tubular neighborhood of the branch, the squared
distance to \(\mathcal M_{t_0}\) is equivalent to the squared norm of the normal
component. We therefore obtain
\[
        \operatorname{dist}((u,v),\mathcal M_{t_0})^2
        \ge c\bigl(\|\xi\|_\alpha^2+\|\zeta\|_\alpha^2+\rho^2\bigr),
\]
provided \(\|\xi\|_\alpha+\|\zeta\|_\alpha+|\rho|\le\varepsilon A\) with
\(\varepsilon\) sufficiently small. The endpoint branches are handled in the
same way. At the branch \((AW_\lambda,0)\), for instance, the tangent
directions are \((W_\lambda,0)\) and \((\partial_\lambda W_\lambda,0)\), while
\((0,W_\lambda)\) is the normal ratio direction.
\end{proof}

\section{Local analysis near the minimizing manifold}
\label{sec:local-analysis}

In this section we prove the sharp local stability estimates near each
minimizing branch. The purpose is to separate the part of the argument which is
common to the non-degenerate and degenerate cases from the part where the two
mechanisms differ.

Let \(t_0\in\mathcal S\cap(0,\infty)\). Every perturbation \((\phi,\psi)\in \mathcal X\) can be written uniquely as
\[
        (\phi,\psi)=e_0\xi+e_1\eta,
\]
or equivalently
\begin{equation}
        \phi=\frac{\xi-t_0\eta}{\sqrt{1+t_0^2}},
        \quad
        \psi=\frac{t_0\xi+\eta}{\sqrt{1+t_0^2}}.
        \label{eq:phi-psi-xi-eta}
\end{equation}
Since \(e_0,e_1\) are orthonormal in \(\mathbb R^2\),
\begin{equation} \label{eq:norm-xi-eta}
        \|(\phi,\psi)\|_{\mathcal X}^2
        =
        \|(\xi,\eta)\|_{\mathcal X}^2 .
\end{equation}
The variable \(\xi\) is the synchronized perturbation, while \(\eta\) is the
transverse perturbation changing the ratio of the two components.

\subsection{Common second-order expansion}\label{5.1}

Throughout Subsection \ref{5.1}, \(t_0\in\mathcal S\cap(0,\infty)\). We first derive a local expansion which will be used in both the non-degenerate
and degenerate cases. By the invariance under \(v\mapsto -v\), it is enough to
consider the positive branch
\[
        (u,v)=A(W_\lambda,t_0W_\lambda)+(\phi,\psi),
        \quad A>0,\quad \lambda>0.
\]

Set
\[
        a_0:=\mu_1+\kappa r t_0^s,
        \quad
        b_0:=\mu_2t_0^{p-1}+\kappa s t_0^{s-1}.
\]
Since \(t_0\) is a critical point of \(g\),
\begin{equation} \label{eq:a0-b0-critical}
        a_0=\frac{\mu_1+\mu_2t_0^p+p\kappa t_0^s}{1+t_0^2},
        \quad
        b_0=t_0a_0.
\end{equation}

Let \(m>0\). Applying Lemmas \ref{lem2.7}, we obtain
\begin{equation}\label{eq:P-common-expansion}
\begin{aligned}
P(AW_\lambda+\phi,At_0W_\lambda+\psi)
\le&
A^p(\mu_1+\mu_2t_0^p+p\kappa t_0^s)
+pA^{p-1}
 \int_{\mathbb R^N}|x|^\ell W_\lambda^{p-1}
       (a_0\phi+b_0\psi)\,dx\\
&+\frac{A^{p-2}}{2}\int_{\mathbb R^N}|x|^\ell W_\lambda^{p-2}
\Big[\big(p(p-1)\mu_1+p\kappa r(r-1)t_0^s
\big)\phi^2\\
&+2p\kappa r s t_0^{s-1}\phi\psi
+\big(p(p-1)\mu_2t_0^{p-2}
    +p\kappa s(s-1)t_0^{s-2}
\big)\psi^2
\Big]\,dx\\
&+\mathcal E_m(A,\phi,\psi),
\end{aligned}
\end{equation}
where
\begin{equation}\label{eq:common-remainder}
\begin{aligned}
\mathcal E_m(A,\phi,\psi)
:=&
mA^{p-2}\int_{\mathbb R^N}|x|^\ell W_\lambda^{p-2}
\Big[
(\mu_1+p\kappa t_0^s)\phi^2
+
(\mu_2t_0^{p-2}+p\kappa t_0^{s-2})\psi^2
\Big]\,dx       \\
&
+C\int_{\mathbb R^N}|x|^\ell\bigl(|\phi|^p+|\psi|^p\bigr)\,dx
+C\int_{\mathbb R^N}|x|^\ell
\mathcal R_{r,s}(AW_\lambda;At_0W_\lambda;\phi,\psi)\,dx .
\end{aligned}
\end{equation}
Moreover,
\begin{equation}\label{eq:energy-common-expansion}
\begin{aligned}
\|AW_\lambda+\phi\|_\alpha^2+\|At_0W_\lambda+\psi\|_\alpha^2
&=
A^2(1+t_0^2)S_0
+2A\langle W_\lambda,\phi\rangle_\alpha
+2At_0\langle W_\lambda,\psi\rangle_\alpha
+\|(\phi,\psi)\|_{\mathcal X}^2 .
\end{aligned}
\end{equation}
The scalar Euler--Lagrange equation gives
\[
        \langle W_\lambda,f\rangle_\alpha
        =
        S_0\int_{\mathbb R^N}|x|^\ell W_\lambda^{p-1}f\,dx,
        \quad f\in X_{\rm rad}.
\]
Hence, by corollary \ref{cor:projection-orthogonality} and \eqref{eq:a0-b0-critical}, the linear terms in
\eqref{eq:P-common-expansion} and \eqref{eq:energy-common-expansion} cancel
exactly. Using
\[
        (a+b)^{2/p}
        \le
        a^{2/p}+\frac2p a^{2/p-1}b,
        \quad a>0,\quad b>-a,
\]
we obtain
\begin{equation}\label{eq:deficit-common-lower-bound}
        \delta(AW_\lambda+\phi,At_0W_\lambda+\psi)
        \ge
        Q_{t_0}(\phi,\psi)
        -C A^{2-p}\mathcal E_m(A,\phi,\psi),
\end{equation}
where
\begin{equation}\label{eq0002}
\begin{aligned}
Q_{t_0}(\phi,\psi)
:=&
\|(\phi,\psi)\|_{\mathcal X}^2 -\frac{(1+t_0^2)S_0}{\mu_1+\mu_2t_0^p+p\kappa t_0^s}
\int_{\mathbb R^N}|x|^\ell W_\lambda^{p-2}
\Big[
\big((p-1)\mu_1+\kappa r(r-1)t_0^s\big)\phi^2 \\
&+2\kappa r s t_0^{s-1}\phi\psi +\big((p-1)\mu_2t_0^{p-2}+\kappa s(s-1)t_0^{s-2}\big)\psi^2
\Big]\,dx .
\end{aligned}
\end{equation}

Substituting
\eqref{eq:phi-psi-xi-eta} into \eqref{eq0002}, and using \(g'(t_0)=0\),
the mixed \(\xi\eta\)-term vanishes. Thus
\begin{equation}\label{eq:Q-diagonalized}
        Q_{t_0}(\phi,\psi)
        =
        Q^{\rm syn}(\xi)
        +
        Q_{t_0}^{\rm tr}(\eta),
\end{equation}
where
\begin{equation}\label{eq:Q-syn}
        Q^{\rm syn}(\xi)
        :=
        \|\xi\|_\alpha^2
        -(p-1)S_0
        \int_{\mathbb R^N}|x|^\ell W_\lambda^{p-2}\xi^2\,dx
\end{equation}
and
\begin{equation}\label{eq:Q-tr}
        Q_{t_0}^{\rm tr}(\eta)
        :=
        \|\eta\|_\alpha^2
        -\Theta(t_0)S_0
        \int_{\mathbb R^N}|x|^\ell W_\lambda^{p-2}\eta^2\,dx,
\end{equation}
where
\[
\Theta(t)
:=
\frac{
(p-1)(\mu_1t^2+\mu_2t^{p-2})
+\kappa\bigl(r(r-1)t^{s+2}-2rs t^s+s(s-1)t^{s-2}\bigr)
}{
\mu_1+\mu_2t^p+p\kappa t^s
}.
\]

We also record the common remainder estimate. Let
\[
        \varrho:=\|(\phi,\psi)\|_{\mathcal X}^2,
        \quad
        \rho_*:=\frac{\varrho^{1/2}}{A}.
\]
If \(\rho_*\to0\), then
\begin{equation}\label{eq:common-remainder-estimate}
        A^{2-p}\mathcal E_m(A,\phi,\psi)
        \le
        Cm\varrho+o(\varrho).
\end{equation}
Indeed, the \(m\)-term is bounded by \(Cm\varrho\). The pure-power remainder satisfies
\[
        A^{2-p}\int_{\mathbb R^N}|x|^\ell
        \bigl(|\phi|^p+|\psi|^p\bigr)\,dx
        \le
        CA^{2-p}\varrho^{p/2}
        =
        C\varrho\rho_*^{p-2}
        =
        o(\varrho).
\]
For the mixed remainder \(\mathcal R_{r,s}\), each term is estimated by the weighted
Sobolev embedding and H\"older's inequality. More precisely,
\[
\begin{array}{ll}
A^{2-p}\displaystyle\int |x|^\ell|\phi|^r|\psi|^s
        \le C\varrho\rho_*^{p-2}=o(\varrho), & r+s=p,\\[1ex]
A^{2-p}\displaystyle\int |x|^\ell(At_0W_\lambda)^s|\phi|^r
        \le C\varrho\rho_*^{r-2}=o(\varrho), & r>2,\\[1ex]
A^{2-p}\displaystyle\int |x|^\ell(AW_\lambda)^r|\psi|^s
        \le C\varrho\rho_*^{s-2}=o(\varrho), & s>2,\\[1ex]
A^{2-p}\displaystyle\int |x|^\ell(AW_\lambda)^{r-2}(At_0W_\lambda)^{s-1}
        |\phi|^2|\psi|
        \le C\varrho\rho_*=o(\varrho), & r\ge2,\\[1ex]
A^{2-p}\displaystyle\int |x|^\ell(AW_\lambda)^{r-1}(At_0W_\lambda)^{s-2}
        |\phi||\psi|^2
        \le C\varrho\rho_*=o(\varrho), & s\ge2,\\[1ex]
        A^{2-p}\displaystyle\int |x|^\ell(AW_\lambda)^{-1}(At_0W_\lambda)^s
        |\phi|^{r+1}
        \le C\varrho\rho_*^{r-1}=o(\varrho), & 1<r<2,\\[1ex]
A^{2-p}\displaystyle\int |x|^\ell(AW_\lambda)^r(At_0W_\lambda)^{-1}
        |\psi|^{s+1}
        \le C\varrho\rho_*^{s-1}=o(\varrho), & 1<s<2,\\[1ex]
A^{2-p}\displaystyle\int |x|^\ell(AW_\lambda)^{-1}(At_0W_\lambda)^{s-1}
        |\phi|^{r+1}|\psi|
        \le C\varrho\rho_*^r=o(\varrho), & 1<r<2\le s,\\[1ex]
A^{2-p}\displaystyle\int |x|^\ell(AW_\lambda)^{r-1}(At_0W_\lambda)^{-1}
        |\phi||\psi|^{s+1}
        \le C\varrho\rho_*^s=o(\varrho), & 1<s<2\le r,\\[1ex]
A^{2-p}\displaystyle\int |x|^\ell(AW_\lambda)^{r-1}(At_0W_\lambda)^{(s-1)/2}
        |\phi||\psi|^{(s+1)/2}
        \le C\varrho\rho_*^{(s-1)/2}=o(\varrho), & 1<r,s<2,\\[1ex]
A^{2-p}\displaystyle\int |x|^\ell(AW_\lambda)^{(r-1)/2}(At_0W_\lambda)^{s-1}
        |\phi|^{(r+1)/2}|\psi|
        \le C\varrho\rho_*^{(r-1)/2}=o(\varrho), & 1<r,s<2.
\end{array}
\]
Thus
\[
        A^{2-p}\int_{\mathbb R^N}|x|^\ell
        \mathcal R_{r,s}(AW_\lambda;At_0W_\lambda;\phi,\psi)\,dx=o(\varrho),
\]
and \eqref{eq:common-remainder-estimate} follows. In particular, for every
\(\varepsilon>0\), by first choosing \(m>0\) sufficiently small and then taking
\(\rho_*>0\) sufficiently small, we have
\begin{equation} \label{eq:common-remainder-eps}
        A^{2-p}\mathcal E_m(A,\phi,\psi)
        \le
        \varepsilon \varrho .
\end{equation}

\subsection{The non-degenerate case}
\begin{proof}[Proof of the non-degenerate part of Theorem \ref{main-thm}]
Assume \(t_0\in \mathcal S\cap(0,\infty)\) is non-degenerate. By Lemma \ref{3.0},
\[
       g'(t_0)=0,\quad g''(t_0)>0 .
\]
A direct differentiation of \(g(t)\) yields
\begin{equation}\label{eq:Theta-gpp}
        1-\Theta(t_0)
        =
        \frac{1+t_0^2}{2}(\mu_1+\mu_2t_0^p+p\kappa t_0^s)^{2/p}g''(t_0).
\end{equation}
Thus \eqref{eq:Theta-gpp} implies
\[
\Theta(t_0)<1.
\]
Hence the transverse quadratic form is coercive:
\begin{equation}\label{eq0003}
        Q_{t_0}^{\rm tr}(\eta)
        \ge
        \min\{1,1-\Theta(t_0)\}\|\eta\|_\alpha^2.
\end{equation}
Indeed, if \(0\le\Theta(t_0)<1\), the first eigenvalue inequality
\(S_0\int |x|^\ell W_\lambda^{p-2}\eta^2\le\|\eta\|_\alpha^2\) gives the
factor \(1-\Theta(t_0)\); if \(\Theta(t_0)<0\), the potential term is
nonnegative and only the factor one is valid. This distinction corrects the
sign-sensitive transverse estimate.
For the synchronized component, by Proposition \ref{prop:scalar-gap-A}, we have
\begin{equation}\label{eq0004}
        Q^{\rm syn}(\xi)
        \ge
        \left(1-\frac{(p-1)S_0}{\Lambda_3}\right)
        \|\xi\|_\alpha^2.
\end{equation}
Combining \eqref{eq0002}-\eqref{eq:Q-tr} and estimates \eqref{eq0003}-\eqref{eq0004}, we obtain
\begin{equation}\label{eq:Q-nondegenerate-coercive}
        Q_{t_0}(\phi,\psi)
        \ge
        c\bigl(\|(\xi,\eta)\|_{\mathcal X}^2\bigr).
\end{equation}
By \eqref{eq:norm-xi-eta}, \eqref{eq:deficit-common-lower-bound},
\eqref{eq:common-remainder-eps}, and \eqref{eq:Q-nondegenerate-coercive},
\[
        \delta(u,v)
        \ge
        c\varrho-\varepsilon \varrho.
\]
Choosing \(\varepsilon>0\) sufficiently small gives
\begin{equation}\label{eq:local-nondegenerate-stability}
        \delta(u,v)
        \ge
        c\varrho\ge c\dist((u,v),{\mathcal M})^2.
\end{equation}

For the endpoint branches, the expansion has to be written with the one-sided ratio.
We spell out the case \(0\in \mathcal S\); the case \(t_0=\infty\) follows by interchanging
the two components. Notice first that \(0\in S\) necessarily implies \(s\ge 2\).
Indeed, if \(1<s<2\), then, as \(t\to0^+\),
\[
g(t)=\mu_1^{-2/p}\left(1+t^2-\frac{2\kappa}{\mu_1}t^s+o(t^s)\right)<g(0),
\]
which contradicts the minimality of \(t=0\). Hence only the cases \(s>2\) and \(s=2\)
can occur at the endpoint \(0\in \mathcal S\).
Similarly, \(\infty\in \mathcal S\) implies \(r\ge 2\), by applying the same argument to
\(\widetilde g(\tau)=g(1/\tau)\) as \(\tau\to0^+\).

Near
\((AW_\lambda,0)\) write
\[
        u=AW_\lambda+\phi,
        \quad
        v=\psi,
\]
with
\[
        \langle \phi,W_\lambda\rangle_\alpha=0,
        \quad
        \langle \phi,\partial_\lambda W_\lambda\rangle_\alpha=0.
\]
The pure \(u\)-part gives the scalar second variation and is controlled by the
radial spectral gap. The mixed term is treated according to the value of
\(s\). If \(s>2\), then
\[
\int |x|^\ell |AW_\lambda+\phi|^r|\psi|^s\,dx
        \le C A^r\|\psi\|_\alpha^s+o(A^{p-2}\|\psi\|_\alpha^2).
\]
After multiplication by \(P(AW_\lambda,0)^{2/p-1}\sim A^{2-p}\), this is bounded by
\[
        C\left(\frac{\|\psi\|_\alpha}{A}\right)^{s-2}\|\psi\|_\alpha^2+o(\|\psi\|_\alpha^2),
\]
and is therefore of higher order in the local regime. If \(s=2\), the mixed term has the
second-order contribution
\[
        p\kappa A^{p-2}
        \int_{\mathbb R^N}|x|^\ell W_\lambda^{p-2}\psi^2\,dx.
\]
Consequently the transverse quadratic form is
\[
        \|\psi\|_\alpha^2
        -\Theta(0)S_0
        \int_{\mathbb R^N}|x|^\ell W_\lambda^{p-2}\psi^2\,dx,
\]
where
\[
        \Theta(0)=
        \begin{cases}
        0, & s>2,\\[1mm]
        \dfrac{2\kappa}{\mu_1}, & s=2.
        \end{cases}
\]
The endpoint non-degeneracy \(g''(0)>0\) is precisely \(\Theta(0)<1\). Hence,
for \(\varrho\ll A^2\),
\begin{equation}\label{eq:01}
        \delta(u,v)
        \ge
        c \varrho\ge c\dist((u,v),{\mathcal M})^2.
\end{equation}
At \(t_0=\infty\) one uses the reflected function
\(\widetilde g(t).\)
\end{proof}

\subsection{The degenerate case}
\begin{proof}[Proof of the degenerate part of Theorem \ref{main-thm}]
We now assume that \(t_0\in \mathcal S\cap(0,\infty)\) is degenerate. By Lemma \ref{3.0},
\begin{equation}\label{eq:degenerate-g}
        g'(t_0)=g''(t_0)=g'''(t_0)=0,
        \quad
        g^{(4)}(t_0)>0 .
\end{equation}
Therefore, identity \eqref{eq:Theta-gpp} gives
\[
        \Theta(t_0)=1.
\]
Thus the transverse quadratic form becomes
\[
        Q_{t_0}^{\rm tr}(\eta)
        =
        \|\eta\|_\alpha^2
        -
        S_0
        \int_{\mathbb R^N}|x|^\ell W_\lambda^{p-2}\eta^2\,dx .
\]
The direction \(W_\lambda\) is exactly the kernel of this form. Therefore the
transverse variable must be decomposed further:
\[
        \eta=\rho W_\lambda+\zeta,
        \quad
        \langle \zeta,W_\lambda\rangle_\alpha=0.
\]
The scalar \(\rho\) is the unique degenerate ratio direction. If
\[
        (\phi,\psi)=e_0\xi+e_1(\rho W_\lambda+\zeta),
\]
with
\begin{equation}\label{eq:02}
        \langle \xi,W_\lambda\rangle_\alpha=
        \langle \xi,\partial_\lambda W_\lambda\rangle_\alpha=
        \langle \zeta,W_\lambda\rangle_\alpha=0,
\end{equation}
then
\begin{equation}\label{eq:Q-degenerate-partial-coercive}
        Q_{t_0}(\phi,\psi)
        \ge
        c\bigl(\|\xi\|_\alpha^2+\|\zeta\|_\alpha^2\bigr).
\end{equation}
Indeed, \(Q^{\rm syn}\) is coercive on
\(\{W_\lambda,\partial_\lambda W_\lambda\}^\perp\), while
\(Q_{t_0}^{\rm tr}\) is coercive on \(W_\lambda^\perp\). The only
direction not controlled by the second variation is
\(
        e_1W_\lambda.
\)

We next control this missing direction by the fourth-order expansion of \(g\).
Consider first the pure degenerate path
\[
        (u_\rho,v_\rho)
        :=
        A(W_\lambda,t_0W_\lambda)+\rho e_1W_\lambda .
\]
Equivalently,
\[
        u_\rho=a_\rho W_\lambda,
        \quad
        v_\rho=b_\rho W_\lambda,
\]
where
\[
        a_\rho=A-\frac{t_0\rho}{\sqrt{1+t_0^2}},
        \quad
        b_\rho=At_0+\frac{\rho}{\sqrt{1+t_0^2}}.
\]
Set
\[
        T_\rho:=\frac{b_\rho}{a_\rho}.
\]
For \(|\rho|\le\varepsilon_0 A\), after reducing \(\varepsilon_0\) if necessary, we have \(a_\rho>0\), \(b_\rho>0\), and \(T_\rho\) remains in a fixed compact subinterval of \((0,\infty)\). Moreover, \(a_\rho\asymp  A\), and
\begin{equation}\label{eq:T-rho-expansion}
        T_\rho-t_0
        =
        \frac{\sqrt{1+t_0^2}\,\rho}
        {A-\dfrac{t_0\rho}{\sqrt{1+t_0^2}}},
        \quad
        |T_\rho-t_0|\asymp \frac{|\rho|}{A}.
\end{equation}
By the definition of \(g\),
\begin{equation}\label{eq:pure-rho-deficit}
\delta(u_\rho,v_\rho)=
S_0a_\rho^2(\mu_1+\mu_2T_{\rho}^p+p\kappa T_{\rho}^s)^{2/p}
\bigl(g(T_\rho)-g(t_0)\bigr).
\end{equation}
Using \eqref{eq:degenerate-g} and Taylor's formula, after shrinking
\(\varepsilon_0\) if necessary, there exists \(c\) such that
\[
        g(T_\rho)-g(t_0)\ge c|T_\rho-t_0|^4.
\]
Therefore, by \eqref{eq:T-rho-expansion},
\begin{equation}\label{eq:pure-rho-control}
        \delta(u_\rho,v_\rho)
        \ge
        cA^2\left(\frac{|\rho|}{A}\right)^4
        =
        c\frac{\rho^4}{A^2}.
\end{equation}

We now add the orthogonal variables \(\xi\) and \(\zeta\):
\[
        (u,v)
        =
        A(W_\lambda,t_0W_\lambda)
        +e_1\rho W_\lambda
        +e_0\xi+e_1\zeta ,
\]
where \(A>0\), \(\lambda>0\), and $\xi,\zeta$ satisfies \eqref{eq:02}.
Equivalently,
\[
        (u,v)=(u_\rho,v_\rho)+(h_1,h_2),
\]
where
\begin{equation}\label{eq:h1-h2-def}
        h_1=\frac{\xi-t_0\zeta}{\sqrt{1+t_0^2}},
        \quad
        h_2=\frac{t_0\xi+\zeta}{\sqrt{1+t_0^2}}.
\end{equation}
By \eqref{eq:02},
\begin{equation}\label{eq:h1-h2-orthogonal}
        \langle h_1,W_\lambda\rangle_\alpha
        =
        \langle h_2,W_\lambda\rangle_\alpha
        =
        0,
\end{equation}
and, since \(e_0,e_1\) are orthonormal,
\begin{equation}\label{eq:h1-h2-norm}
        \|(h_1,h_2)\|_{\mathcal X}^2
        =
        \|(\xi,\zeta)\|_{\mathcal X}^2.
\end{equation}
Throughout the following argument we assume
\begin{equation}
        \|\xi\|_\alpha+\|\zeta\|_\alpha+|\rho|
        \le
        \varepsilon_0 A,
        \label{eq:degenerate-smallness}
\end{equation}
where \(\varepsilon_0>0\) will be chosen sufficiently small. Then
\[
        a_\rho\sim A,
        \quad
        b_\rho\sim A,
        \quad
        P(u_\rho,v_\rho)\sim A^p.
\]

We start from the exact decomposition of the deficit. Since the energy is
quadratic and the linear energy terms vanish by \eqref{eq:h1-h2-orthogonal},
\[
\begin{aligned}
\delta(u,v)
=&
\delta(u_\rho,v_\rho)
+\|(h_1,h_2)\|_{\mathcal X}^2 \\
&-S_{r,s,\mu_1,\mu_2,\kappa}
\left[
        P(u_\rho+h_1,v_\rho+h_2)^{2/p}
        -
        P(u_\rho,v_\rho)^{2/p}
\right].
\end{aligned}
\]
By the concavity of \(t\mapsto t^{2/p}\),
\begin{equation}\label{eq:deficit-concavity-degenerate}
\begin{aligned}
\delta(u,v)
\ge&
\delta(u_\rho,v_\rho)
+\|(h_1,h_2)\|_{\mathcal X}^2 \\
&
-\frac{2S_{r,s,\mu_1,\mu_2,\kappa}}{p}
P(u_\rho,v_\rho)^{2/p-1}
\left[
        P(u_\rho+h_1,v_\rho+h_2)-P(u_\rho,v_\rho)
\right].
\end{aligned}
\end{equation}

We next expand the nonlinear increment in \eqref{eq:deficit-concavity-degenerate}
up to second order. Lemma \ref{lem2.7} and \eqref{eq:h1-h2-orthogonal} gives, for
any \(m>0\),
\begin{equation}\label{eq:pure-first-degenerate}
\begin{aligned}
&P(u_\rho+h_1,v_\rho+h_2)-P(u_\rho,v_\rho)\\
\leq&\left(\frac{p(p-1)}2+m\right)
\int_{\mathbb R^N}|x|^\ell W_\lambda^{p-2}(\mu_1|a_\rho|^{p-2}h_1^2+\mu_2|b_\rho|^{p-2}h_2^2)\,dx \\
&+C_m\int_{\mathbb R^N}|x|^\ell (|h_1|^p+|h_2|^p)\,dx\\
&+p\kappa\int_{\mathbb R^N}|x|^\ell \left(\frac{r(r-1)}2+m\right)
|a_\rho|^{r-2}|b_\rho|^s
W_\lambda^{p-2}h_1^2\,dx  \\
&+p\kappa\int_{\mathbb R^N}|x|^\ell\left(\frac{s(s-1)}2+m\right)
|a_\rho|^r|b_\rho|^{s-2}
W_\lambda^{p-2}h_2^2\,dx \\
&+p\kappa rs\int_{\mathbb R^N}|x|^\ell|a_\rho|^{r-2}a_\rho |b_\rho|^{s-2}b_\rho
W_\lambda^{p-2}h_1h_2\,dx \\
&+
C_m\int_{\mathbb R^N}|x|^\ell  \mathcal  R_{r,s}(a_\rho W_\lambda,b_\rho W_\lambda;h_1,h_2)\,dx.
\end{aligned}
\end{equation}

The remainders in the pure terms are higher order after multiplication by
\(P(u_\rho,v_\rho)^{2/p-1}\). Indeed, by the scalar weighted Sobolev inequality,
\begin{equation}\label{eq:pure-high-order-small}
P(u_\rho,v_\rho)^{2/p-1}
\int_{\mathbb R^N}|x|^\ell |h_i|^p\,dx
\le
C\left(\frac{\|h_i\|_\alpha}{A}\right)^{p-2}\|h_i\|_\alpha^2
=
o(1)\|h_i\|_\alpha^2,
\quad i=1,2.
\end{equation}
The \(m\)-terms are quadratic and satisfy
\begin{align}
&P(u_\rho,v_\rho)^{2/p-1}
m\mu_1|a_\rho|^{p-2}
\int_{\mathbb R^N}|x|^\ell W_\lambda^{p-2}h_1^2\,dx
\le
Cm\|h_1\|_\alpha^2,
\label{eq:mixed-m-term-1}
\\
&P(u_\rho,v_\rho)^{2/p-1}
m\mu_2|b_\rho|^{p-2}
\int_{\mathbb R^N}|x|^\ell W_\lambda^{p-2}h_2^2\,dx
\le
Cm\|h_2\|_\alpha^2,
\label{eq:mixed-m-term-2}
\\
&P(u_\rho,v_\rho)^{2/p-1}
m|a_\rho|^{r-2}|b_\rho|^s
\int_{\mathbb R^N}|x|^\ell W_\lambda^{p-2}h_1^2\,dx
\le
Cm\|h_1\|_\alpha^2,
\label{eq:mixed-m-term-3}
\\
&P(u_\rho,v_\rho)^{2/p-1}
m|a_\rho|^r|b_\rho|^{s-2}
\int_{\mathbb R^N}|x|^\ell W_\lambda^{p-2}h_2^2\,dx
\le
Cm\|h_2\|_\alpha^2.
\label{eq:mixed-m-term-4}
\end{align}
As in non-degenerate case, the higher-order terms satisfy
\begin{equation}\label{eq:mixed-high-order-small}
P(u_\rho,v_\rho)^{2/p-1}\mathcal R_{r,s}(a_\rho W_\lambda,b_\rho W_\lambda;h_1,h_2)=o(1)\bigl(\|(h_1,h_2)\|_{\mathcal X}^2\bigr).
\end{equation}
Combining \eqref{eq:deficit-concavity-degenerate}--\eqref{eq:mixed-high-order-small},
we obtain
\begin{equation}\label{eq:delta-after-expansion}
\delta(u,v)
\ge
\delta(u_\rho,v_\rho)
+Q_\rho(h_1,h_2)
-Cm\bigl(\|(h_1,h_2)\|_{\mathcal X}^2\bigr)
-o(1)\bigl(\|(h_1,h_2)\|_{\mathcal X}^2\bigr),
\end{equation}
where
\begin{equation}\label{eq:explicit-Qrho}
\begin{aligned}
Q_\rho(h_1,h_2)
:=&
\|(h_1,h_2)\|_{\mathcal X}^2
-
S_{r,s,\mu_1,\mu_2,\kappa}(\mu_1+\mu_2T_{\rho}^p+p\kappa T_{\rho}^s)^{2/p-1}\\
&\times
\int_{\mathbb R^N}|x|^\ell W_\lambda^{p-2}
\big[
((p-1)\mu_1+\kappa r(r-1)T_{\rho}^s)h_1^2\\
&
+2\kappa rs\,T_{\rho}^{s-1}h_1h_2
+((p-1)\mu_2T_{\rho}^{p-2}+\kappa s(s-1)T_{\rho}^{s-2})h_2^2
\big]\,dx,
\end{aligned}
\end{equation}

Since \(\frac{|\rho|}{A}
\) is  sufficiently small
and the coefficients of the above quadratic form depend smoothly on the
ratio \(T_\rho\) for \(T_\rho\) close to \(t_0\), we have
\begin{equation}\label{eq:Qr-Q0-comparison}
        \left|
        Q_\rho(h_1,h_2)
        -
        Q_{t_0}(h_1,h_2)
        \right|
        \le
        C\frac{|\rho|}{A}
        \bigl(\|(h_1,h_2)\|_{\mathcal X}^2\bigr)
        =
        o(1)\bigl(\|(h_1,h_2)\|_{\mathcal X}^2\bigr).
\end{equation}
Therefore,
\begin{equation}\label{eq:delta-with-Q-t0}
\delta(u,v)
\ge
\delta(u_\rho,v_\rho)
+Q_{t_0}(h_1,h_2)
-Cm\bigl(\|(h_1,h_2)\|_{\mathcal X}^2\bigr)
-o(1)\bigl(\|(h_1,h_2)\|_{\mathcal X}^2\bigr).
\end{equation}

By the partial coercivity estimate \eqref{eq:Q-degenerate-partial-coercive}, applied with \(\rho=0\), we have
\begin{equation}\label{eq:03}
        Q_{t_0}(h_1,h_2)
        =
        Q_{t_0}(e_0\xi+e_1\zeta)
        \ge
        c\bigl(\|\xi\|_\alpha^2+\|\zeta\|_\alpha^2\bigr).
\end{equation}
Combining \eqref{eq:pure-rho-control}, \eqref{eq:h1-h2-norm} and \eqref{eq:03}, choosing \(m>0\) sufficiently small and
then taking \(\varepsilon_0>0\) sufficiently small, all error terms in
\eqref{eq:delta-with-Q-t0} can be absorbed. Hence
\begin{equation} \label{eq:local-degenerate-stability-final}
        \delta(u,v)
        \ge
        c\left(
        \|\xi\|_\alpha^2+\|\zeta\|_\alpha^2+\frac{\rho^4}{A^2}
        \right).
\end{equation}

By the local smallness assumption and the decomposition around the branch, we have
\[\|(u,v)\|_{\mathcal X}^2=A^2(1+t_0^2)S_0+o(A^2),\]
and hence
\begin{equation} \label{eq:001}
        A^2\asymp \|(u,v)\|_{\mathcal X}^2.
\end{equation}
By Lemma \ref{lem:local-distance}, we have
\begin{equation}\label{eq:002}
        \operatorname{dist}\bigl((u,v),\mathcal M_{t_0}\bigr)^2
        \asymp
        \|\xi\|_\alpha^2+\|\zeta\|_\alpha^2+\rho^2,
\end{equation}
Combining \eqref{eq:local-degenerate-stability-final},  \eqref{eq:001} and \eqref{eq:002},
we obtain
\[
        \frac{\delta(u,v)}{\|(u,v)\|_{\mathcal X}^2}
        \ge
        c\left[
        \frac{\|\xi\|_\alpha^2+\|\zeta\|_\alpha^2}{\|(u,v)\|_{\mathcal X}^2}
        +
        \left(\frac{\rho^2}{\|(u,v)\|_{\mathcal X}^2}\right)^2
        \right]\geq c\left(
        \frac{\operatorname{dist}\bigl((u,v),\mathcal M_{t_0}\bigr)^2}{\|(u,v)\|_{\mathcal X}^2}
        \right)^2.
\]
Therefore,
\[
        \left(
        \frac{\delta(u,v)}
        {\|(u,v)\|_{\mathcal X}^2}
        \right)^{1/2}
        \ge
        C
        \frac{\operatorname{dist}\bigl((u,v),\mathcal M_{t_0}\bigr)^2}
        {\|(u,v)\|_{\mathcal X}^2}
        \ge
        C \frac{\operatorname{dist}\bigl((u,v),\mathcal M\bigr)^2}
        {\|(u,v)\|_{\mathcal X}^2}.
\]

If \(t_0=0\), the reference bubble is \((AW_\lambda,0)\), and we write
\[
        u=AW_\lambda+\xi,
        \quad
        v=\rho W_\lambda+\zeta,
\]
with
\[
        \langle \xi,W_\lambda\rangle_\alpha=0,
        \quad
        \langle \xi,\partial_\lambda W_\lambda\rangle_\alpha=0,
        \quad
        \langle \zeta,W_\lambda\rangle_\alpha=0.
\]
The pure degenerate path \((AW_\lambda,\rho W_\lambda)\) corresponds to the
one-sided ratio \(t=|\rho|/A\), since the functional depends on \(|v|\). If \(0\in \mathcal S\) is degenerate, then
\[
        g(t)-g(0)\ge ct^4
        \quad
        \text{for }0\le t\ll1.
\]
Consequently,
\[
        \delta(AW_\lambda,\rho W_\lambda)
        \ge
        cA^2\left(\frac{|\rho|}{A}\right)^4
        =
        c\frac{\rho^4}{A^2}.
\]

The orthogonal variables \(\xi\) and \(\zeta\) are controlled by the scalar
spectral gap exactly as above. The case \(t_0=\infty\) is obtained by
interchanging the two components. In this case the pure degenerate path is
\((\rho W_\lambda,AW_\lambda)\), and the one-sided ratio is \(t=|\rho|/A\) for
\(
\widetilde g(t).
\)
\end{proof}

\begin{remark}\label{rem:sharpness-exponent}
The exponents \(\iota=1\) and \(\iota=1/2\) in Theorem~\ref{main-thm} for non-flat case are sharp.
Let
\(
(u_t,v_t):=(W_{\lambda},t W_{\lambda}).
\)
If \(t_0\in \mathcal S\cap[0,\infty)\), then
\[
\frac{\operatorname{dist}((u_t,v_t),\mathcal M)^2}
{\|u_t\|_\alpha^2+\|v_t\|_\alpha^2}
\asymp |t-t_0|^2
\quad\text{as }t\to t_0.
\]
Moreover,
\[
\delta(u_t,v_t)
=
S_{0}(\mu_1+\mu_2t^p+p\kappa t^s)^{2/p}\bigl(g(t)-m_*\bigr).
\]
Thus,
\[
\frac{\delta(u_t,v_t)}
{\|u_t\|_\alpha^2+\|v_t\|_\alpha^2}
\asymp |t-t_0|^2,
\quad\text{as }t\to t_0
\]
in the non-degenerate case, and
\[
\frac{\delta(u_t,v_t)}
{\|u_t\|_\alpha^2+\|v_t\|_\alpha^2}
\asymp |t-t_0|^4,
\quad\text{as }t\to t_0
\]
in the degenerate case.

If \(\infty\in \mathcal S\), one instead takes
\[
(u_t,v_t)=(tW_{\lambda},W_{\lambda}),\quad t\rightarrow0^{+},
\]
and uses the function \(\widetilde g(t)=g(1/t)\).
\end{remark}

\subsection{The flat case}\label{sec:flat}
Since \( (r,s,\mu_1,\mu_2)=(2,2,2\kappa,2\kappa)\),
\begin{equation}\label{flat-deficit}
\delta(u,v)
:=
\|(u,v)\|_{\mathcal X}^2
-
S_0
\left(
\int_{\mathbb R^N}|x|^\ell
\bigl(|u|^2+|v|^2\bigr)^2\,dx
\right)^{1/2},
\end{equation}
The functional \(\delta\) is invariant under the natural action of
\(O(2)\): if \(R\in O(2)\) is constant and
\[
\binom{\widetilde u}{\widetilde v}
=
R\binom{u}{v},
\]
then
\[
\|(\widetilde u,\widetilde v)\|_{\mathcal X}
=
\|(u,v)\|_{\mathcal X},
\quad
\delta(\widetilde u,\widetilde v)=\delta(u,v).
\]
This additional symmetry is responsible for the flat ratio direction.

First, we treat the large-deficit regime. Since
\[
\operatorname{dist}((u,v),\mathcal M_{\mathrm{flat}})^2
\leq \|(u,v)\|_X^2,
\]
if
\[
\delta(u,v)\geq \varepsilon_0\|(u,v)\|_X^2,
\]
then
\[
\delta(u,v)
\geq
\varepsilon_0
\operatorname{dist}((u,v),\mathcal M_{\mathrm{flat}})^2 .
\]
Therefore, after changing the stability constant if necessary, it remains to consider the small-deficit regime.

It is convenient to write
\[
\mathcal M_{\mathrm{flat}}
=
\left\{
A e W_\lambda:
A>0,\ e\in\mathbb S^1,\ \lambda>0
\right\},
\]
where \(e=(e^1,e^2)\in\mathbb S^1\) and
\[
eW_\lambda=(e^1W_\lambda,e^2W_\lambda).
\]

\begin{lemma}
\label{lem:flat-projection}
Suppose that
\[
\operatorname{dist}
\bigl((u,v),\mathcal M_{\mathrm{flat}}\bigr)
<
\|(u,v)\|_{\mathcal X}.
\]
Then the distance to \(\mathcal M_{\mathrm{flat}}\) is attained at some
point
\[
A e W_\lambda\in\mathcal M_{\mathrm{flat}},
\qquad
A>0,\quad e\in\mathbb S^1,\quad\lambda>0.
\]
After applying a constant orthogonal transformation in the component
variables, one may assume that \(e=(1,0)\), so that
\[
u=AW_\lambda+\phi,
\quad
v=\psi,
\]
where
\begin{equation}\label{eq-01}
\langle\phi,W_\lambda\rangle_\alpha=0,
\quad
\langle\phi,\partial_\lambda W_\lambda\rangle_\alpha=0,
\quad
\langle\psi,W_\lambda\rangle_\alpha=0.
\end{equation}
Moreover,
\begin{equation}\label{eq-02}
\operatorname{dist}
\bigl((u,v),\mathcal M_{\mathrm{flat}}\bigr)^2
=
\|(\phi,\psi)\|_{\mathcal X}^2.
\end{equation}
\end{lemma}

\begin{proof}
The existence of a minimizing projection in a sufficiently small
neighborhood follows from the same weak-vanishing argument as in
Lemma~\ref{lem:projection-attainment}, together with compactness of the angular variable
\(e\in\mathbb S^1\).

Since the functional, the norm, and the minimizing manifold are
invariant under \(O(2)\), we may rotate the minimizing direction \(e\)
to \((1,0)\). The tangent space to
\(\mathcal M_{\mathrm{flat}}\) at \((AW_\lambda,0)\) is
\[
T_{(AW_\lambda,0)}\mathcal M_{\mathrm{flat}}
=
\operatorname{span}
\left\{
(W_\lambda,0),
(\partial_\lambda W_\lambda,0),
(0,W_\lambda)
\right\}.
\]
Differentiating the squared distance with respect to the amplitude,
scale, and angular parameters gives precisely
\eqref{eq-01}. Since the chosen point is a nearest point,
\eqref{eq-02} follows.
\end{proof}

We next establish a qualitative compactness result.
\begin{proposition}\label{prop:flat-compactness}
Let \((u_n,v_n)\subset\mathcal X\) satisfy
\(
\|(u_n,v_n)\|_{\mathcal X}^2=1
\)
and
\(
\delta(u_n,v_n)\longrightarrow0.
\)
Then
\begin{equation}\label{eq-03}
\operatorname{dist}
\bigl((u_n,v_n),\mathcal M_{\mathrm{flat}}\bigr)
\longrightarrow0.
\end{equation}
\end{proposition}

\begin{proof}
The proof follows the same compactness argument as in Proposition \ref{prop:system-compactness}.
We only indicate the modifications in the present flat case.
Define $a_n,b_n$ as in the proof of Proposition \ref{prop:system-compactness}.
Since
\(
\delta(u_n,v_n)\longrightarrow 0,
\)
all inequalities used in the proof of Proposition \ref{prop:system-compactness} become asymptotically equalities.
In particular,
\[
\|u_n\|_\alpha^2-S_0a_n^2\longrightarrow 0,
\quad
\|v_n\|_\alpha^2-S_0b_n^2\longrightarrow 0.
\]

By the scalar compactness of minimizing sequences, after passing to a subsequence,
\[
u_n=c_nW_{\lambda_n}+o(1),
\quad
v_n=d_nW_{\mu_n}+o(1),\quad \text{in} \quad X_{\mathrm{rad}}.
\]
The equality case in H\"older's inequality yields
\[
\frac{\lambda_n}{\mu_n}\longrightarrow 1,
\]
and hence
\[
\operatorname{dist}\big((u_n,v_n),\mathcal M_{\mathrm{flat}}\big)
\longrightarrow 0.
\]
\end{proof}

We first establish local coercivity.
\begin{proposition}
\label{prop:flat-local-stability}
There exist constants \(\varepsilon_0>0\) and \(c>0\) such that, if
\(
(u,v)=(AW_\lambda+\phi,\ \psi)
\)
with the orthogonality conditions \eqref{eq-01}, and
\(
\|(\phi,\psi)\|_{\mathcal X}
\leq \varepsilon_0 A,
\)
then
\begin{equation}\label{eq-05}
\delta(u,v)
\geq
c
\left(
\|(\phi,\psi)\|_{\mathcal X}^2
\right).
\end{equation}
\end{proposition}

\begin{proof}
A direct computation shows that
\begin{equation}\label{eq-06}
\begin{aligned}
\delta(u,v)=&\|AW+\phi\|_\alpha^2+\|\psi\|_\alpha^2-S_0\left(
\int_{\mathbb R^N}|x|^\ell
\bigl(A^2W^2+2AW\phi+\phi^2+\psi^2\bigr)^2\,dx
\right)^{1/2}\\
=&\|(\phi,\psi)\|_{\mathcal X}^2-3S_0\int_{\mathbb R^N}|x|^\ell W^2\phi^2\,dx
-S_0\int_{\mathbb R^N}|x|^\ell W^2\psi^2\,dx-\mathcal R(A,\phi,\psi),
\end{aligned}
\end{equation}
where
\begin{equation}\label{eq-07}
|\mathcal R(A,\phi,\psi)|
\leq
C
\left[
\frac{
\bigl(\|(\phi,\psi)\|_{\mathcal X}\bigr)^3
}{A}
+
\frac{
\bigl(\|(\phi,\psi)\|_{\mathcal X}\bigr)^4
}{A^2}
\right]=o(\|(\phi,\psi)\|_{\mathcal X}^2),
\end{equation}
as \(\varepsilon_0\) is sufficiently small.
By Proposition~\ref{prop:scalar-gap-A}, \eqref{eq-06} and \eqref{eq-07}, we have
\[
\delta(u,v)\geq\left(
1-\frac{3S_0}{\Lambda_3}
\right)
\|\phi\|_\alpha^2+\frac23\|\psi\|_\alpha^2
+o(\|(\phi,\psi)\|_{\mathcal X}^2)\geq c\|(\phi,\psi)\|_{\mathcal X}^2.
\]
\end{proof}

\begin{proof}[Proof of Theorem~\ref{main-thm} in the flat case]
By homogeneity, it is enough to prove the estimate under the
normalization
\(
\|(u,v)\|_{\mathcal X}^2=1.
\)

Suppose, by contradiction, that the conclusion is false. Then there
exists a sequence \((u_n,v_n)\subset\mathcal X\) such that
\[
\|(u_n,v_n)\|_{\mathcal X}^2=1
\]
and
\begin{equation}\label{eq-08}
\frac{
\delta(u_n,v_n)
}{
\operatorname{dist}
\bigl((u_n,v_n),\mathcal M_{\mathrm{flat}}\bigr)^2
}
\longrightarrow0.
\end{equation}
Since
\[
\operatorname{dist}
\bigl((u_n,v_n),\mathcal M_{\mathrm{flat}}\bigr)^2
\leq
\|(u_n,v_n)\|_{\mathcal X}^2=1,
\]
relation \eqref{eq-08} implies
\[
\delta(u_n,v_n)\longrightarrow0.
\]
Proposition~\ref{prop:flat-compactness} therefore gives
\[
\operatorname{dist}
\bigl((u_n,v_n),\mathcal M_{\mathrm{flat}}\bigr)
\longrightarrow0.
\]
For all sufficiently large \(n\), Lemma~\ref{lem:flat-projection}
allows us, after a constant \(O(2)\)-rotation, to write
\[
u_n=A_nW_{\lambda_n}+\phi_n,
\qquad
v_n=\psi_n,
\]
with
\[
\langle\phi_n,W_{\lambda_n}\rangle_\alpha
=
\langle\phi_n,\partial_\lambda W_{\lambda_n}\rangle_\alpha
=
\langle\psi_n,W_{\lambda_n}\rangle_\alpha
=0
\]
and
\[
\|\phi_n\|_\alpha^2+\|\psi_n\|_\alpha^2
=
\operatorname{dist}
\bigl((u_n,v_n),\mathcal M_{\mathrm{flat}}\bigr)^2.
\]
The normalization and the convergence to the minimizing manifold imply
that \(A_n\) is bounded away from zero and
\[
\frac{\|\phi_n\|_\alpha+\|\psi_n\|_\alpha}{A_n}
\longrightarrow0.
\]
Proposition~\ref{prop:flat-local-stability} then yields
\[
\delta(u_n,v_n)
\geq
c_0
\operatorname{dist}
\bigl((u_n,v_n),\mathcal M_{\mathrm{flat}}\bigr)^2,
\]
contradicting \eqref{eq-08}.

It remains to prove the optimal of the stability exponent statement.
Let
\((u_\varepsilon,v_\varepsilon):=
(aW_{\lambda}+\varepsilon\phi,bW_{\lambda})
\) with \((a,b)\neq(0,0)\),
where \(\phi\) satisfy
\begin{equation}\label{eq-007}
\langle\phi,W_\lambda\rangle_\alpha=0,
\quad
\langle\phi,\partial_\lambda W_\lambda\rangle_\alpha=0, \quad \|\phi\|_{\alpha}=1.
\end{equation}
By \eqref{flat-deficit} and \eqref{eq-007}, we have
\[
\begin{aligned}
\delta(u_\varepsilon,v_\varepsilon)=&
(a^2+b^2)S_0+\varepsilon^2\|\phi\|_\alpha^2-S_0\left(
\int_{\mathbb R^N}
|x|^\ell
\bigl(|u_\varepsilon|^2+|v_\varepsilon|^2\bigr)^2\,dx
\right)^{1/2}\\
=&\varepsilon^2
\left[
\|\phi\|_\alpha^2
-
S_0\frac{3a^2+b^2}{a^2+b^2}
\int_{\mathbb R^N}|x|^\ell W_1^2\phi^2\,dx
\right]
+
O(|\varepsilon|^3).
\end{aligned}
\]
By Proposition \ref{prop:scalar-gap-A} and \ref{eq-007}, we have
\[
\|\phi\|_\alpha^2
\geq
\Lambda_3\int_{\mathbb R^N}|x|^\ell W_1^2\phi^2\,dx>S_0\frac{3a^2+b^2}{a^2+b^2}
\int_{\mathbb R^N}|x|^\ell W_1^2\phi^2\,dx.
\]
 It follows that
\[
\delta(u_\varepsilon,v_\varepsilon)\asymp\varepsilon^2.
\]
On the other hand,
\[
\operatorname{dist}
\bigl((u_\varepsilon,v_\varepsilon),\mathcal M_{\rm flat}\bigr)^2
\asymp\varepsilon^2.
\]
Thus the exponent in the flat case of Theorem~\ref{main-thm} is optimal.
\end{proof}

\medskip
\noindent\textbf{Data availability.}
No data sets were generated or analyzed in this study.

\noindent\textbf{Conflict of interest.}
The authors declare that they have no conflict of interest.


\begin{thebibliography}{99}

\bibitem{BE1991}
G. Bianchi, H. Egnell,
\textit{A note on the Sobolev inequality},
J. Funct. Anal., {\bf 100} (1) (1991), 18--24.

\bibitem{BL1985}
H. Brezis, E. H. Lieb,
\textit{Sobolev inequalities with remainder terms},
J. Funct. Anal., {\bf 62} (1) (1985), 73--86.

\bibitem{CKN1984}
L. Caffarelli, R. Kohn, L. Nirenberg,
\textit{First order interpolation inequalities with weights},
Compos. Math., {\bf 53} (3) (1984), 259--275.

\bibitem{CS2025}
S. Chakraborty, U. Sarkar,
\textit{Quantitative stability in fractional Hardy--Sobolev inequalities: the role of Euler--Lagrange equations},
Proc. Roy. Soc. Edinburgh Sect. A, published online first, 2025, 1--39.


\bibitem{CFL2024}
C. Cazacu, J. Flynn, N. Lam, G. Lu,
\textit{Caffarelli--Kohn--Nirenberg identities, inequalities and their stabilities},
J. Math. Pures Appl., {\bf 182} (2024), 253--284.

\bibitem{CKW2024}
H. Chen, S. Kim, J. Wei,
\textit{Sharp quantitative stability estimates for critical points of fractional Sobolev inequalities},
Int. Math. Res. Not. IMRN, {\bf 2025} (12) (2025), rnaf156.

\bibitem{DGK2025}
N. De Nitti, F. Glaudo, T. K\"onig,
\textit{Non-degeneracy, stability, and symmetry for the fractional Caffarelli--Kohn--Nirenberg inequality},
Comm. Partial Differential Equations, {\bf 50} (8) (2025), 984--1028.

\bibitem{DGT2023}
S. Deng, M. Grossi, X. Tian,
\textit{On some weighted fourth order equations},
J. Differential Equations, {\bf364} (2023), 612--634.

\bibitem{DT2024}
S. Deng, X. Tian,
\textit{Classification and non-degeneracy of positive radial solutions for a weighted fourth-order equation and its application},
Nonlinear Anal., {\bf 240} (2024), 113468.

\bibitem{DT2025}
S. Deng, X. Tian,
\textit{Stability of Rellich--Sobolev type inequality involving Hardy term for bi-Laplacian},
J. Differential Equations, {\bf 422} (2025), 78--105.

\bibitem{DT2023}
S. Deng, X. Tian,
\textit{Some weighted fourth-order Hardy--H\'enon equations},
J. Funct. Anal., {\bf 284} (2023), 109745.

\bibitem{FN2019}
A. Figalli, R. Neumayer,
\textit{Gradient stability for the Sobolev inequality: the case $p\geq 2$},
J. Eur. Math. Soc., {\bf 21} (2) (2019), 319--354.

\bibitem{FZ2022}
A. Figalli, Y. R.-Y. Zhang,
\textit{Sharp gradient stability for the Sobolev inequality},
Duke Math. J., {\bf 171} (12) (2022), 2407--2459.

\bibitem{FS2003}
V. Felli, M. Schneider,
\textit{Perturbation results of critical elliptic equations of Caffarelli--Kohn--Nirenberg type},
J. Differential Equations, {\bf 191} (1) (2003), 121--142.

\bibitem{F2022}
R. L. Frank,
\textit{Degenerate stability of some Sobolev inequalities},
Ann. Inst. H. Poincar\'e C Anal. Non Lin\'eaire, {\bf 39} (6) (2022), 1459--1484.

\bibitem{FP2024}
R. L. Frank, J. W. Peteranderl,
\textit{Degenerate stability of the Caffarelli--Kohn--Nirenberg inequality along the Felli--Schneider curve},
Calc. Var. Partial Differential Equations, {\bf 63} (44) (2024).

\bibitem{K2023}
T. K\"onig,
\textit{On the sharp constant in the Bianchi--Egnell stability inequality},
Bull. Lond. Math. Soc., {\bf 55} (4) (2023), 2070--2075.

\bibitem{K2025}
T. K\"onig,
\textit{Stability for the Sobolev inequality: existence of a minimizer},
J. Eur. Math. Soc., published online first, (2025).

\bibitem{LZ2025}
G. Liu, Y. Ru, Y. Zhang,
\textit{Sharp stability for critical points of the Sobolev inequality in the absence of bubbling},
preprint, arXiv:2503.02340, 2025.


\bibitem{LW2000}
G. Lu, J. Wei,
\textit{On a Sobolev inequality with remainder terms},
Proc. Amer. Math. Soc., {\bf 128} (1) (2000), 75--84.


\bibitem{ZW2024}
Y. Zhou, W. Zou,
\textit{Degenerate stability of critical points of the Caffarelli--Kohn--Nirenberg inequality along the Felli--Schneider curve},
J. Lond. Math. Soc. (2), {\bf 113} (2) (2026), e70454.

\bibitem{WW2022}
J. Wei, Y. Wu,
\textit{On the stability of the Caffarelli--Kohn--Nirenberg inequality},
Math. Ann., {\bf 384} (3-4) (2022), 1509--1546.

\bibitem{WW2024}
J. Wei, Y. Wu,
\textit{Stability of the Caffarelli--Kohn--Nirenberg inequality: the existence of minimizers},
Math. Z., {\bf 308} (4) (2024), 64.


\bibitem{WW2026}
J. Wei, Y. Wu,
\textit{Stability of the Caffarelli--Kohn--Nirenberg inequality along the Felli--Schneider curve: critical points at infinity},
Proc. Lond. Math. Soc., {\bf 132} (4) (2026), e70146.

\bibitem{ZZZ2026}
Y. Zhang, X. Zhong, W. Zou,
\textit{Bivariate Hardy--Sobolev inequality and its sharp stability},
preprint, arXiv:2602.03191, 2026.

\bibitem{ZZ2026}
Y. Zhou, W. Zou,
\textit{The stability on the Caffarelli--Kohn--Nirenberg and Hardy-type inequalities and beyond},
J. Differential Equations, {\bf 450} (2026), 113738.





\end{thebibliography}
\end{document}